\documentclass[preprint,11pt]{elsarticle}
\biboptions{sort&compress}
\usepackage[a4paper,margin=30mm]{geometry}
\usepackage[T1]{fontenc}
\usepackage{lmodern}
\usepackage{amsmath,amssymb,amsthm,mathtools}
\usepackage{booktabs,array,graphicx,tikz,flafter}
\graphicspath{{./}{figures/}}
\usepackage{microtype,float}
\usepackage[hidelinks]{hyperref}
\numberwithin{equation}{section}
\newtheorem{theorem}{Theorem}[section]
\newtheorem{proposition}[theorem]{Proposition}
\newtheorem{lemma}[theorem]{Lemma}
\newtheorem{corollary}[theorem]{Corollary}
\theoremstyle{definition}

\newtheorem{remark}[theorem]{Remark}
\newcommand{\C}{\mathbb C}
\newcommand{\R}{\mathbb R}

\newcommand{\Chat}{\widehat{\mathbb C}}

\newcommand{\RR}{\mathcal R}
\newcommand{\Kell}{\mathcal K}
\newcommand{\Eell}{\mathcal E}
\newcommand{\supp}{\operatorname{supp}}
\newcommand{\capc}{\operatorname{cap}}
\newcommand{\dist}{\operatorname{dist}}
\newcommand{\weakto}{\stackrel{*}{\rightharpoonup}}
\newcommand{\eps}{\varepsilon}

\newcommand{\Rea}{\operatorname{Re}}
\newcommand{\Ima}{\operatorname{Im}}
\journal{Journal of Mathematical Analysis and Applications}
\hypersetup{pdftitle={Pole rigidity in root-optimal rational approximation of the exponential on infinite sectors},
 pdfauthor={Fei Xue and Tianqi Zhang},
 pdfkeywords={rational approximation, pole distribution, infinite sector, Green equilibrium, S-property, common denominator}}
\begin{document}
\begin{frontmatter}
\title{Pole rigidity in root-optimal rational approximation of the exponential on infinite sectors}
\author[clemson]{Fei Xue\corref{cor}}
\ead{fxue@clemson.edu}
\author[zufe]{Tianqi Zhang}
\ead{tianqiz@zufe.edu.cn}
\cortext[cor]{Corresponding author.}
\address[clemson]{School of Mathematical and Statistical Sciences, Clemson University, 220 Parkway Dr., Clemson, SC 29634, USA}
\address[zufe]{School of Data Sciences, Zhejiang University of Finance and Economics, 18 Xueyuan Street, Hangzhou, Zhejiang 310018, China}
\begin{abstract}
The optimal root convergence factor for uniform rational approximation of the exponential on an infinite sector is classical. We prove that, for every fixed infinite sector of opening less than $\pi$, the normalized pole-counting measures of any reduced root-optimal sequence, after scaling the poles by degree, converge weakly to a uniquely determined probability measure supported on a compact analytic arc outside the sector. For the half-line, the corresponding pole law was previously announced by Stahl and Schmelzer; our new result covers sectors of positive opening. The proof rules out boundary pole accumulation and controls an exponentially small Hermite residual through reflection and noncancellation, without assuming interpolation or orthogonality. We explicitly construct the equilibrium measure and its supporting arc in terms of a zero of a Ferrers function, and relate this parameterization to the Gonchar--Rakhmanov theta equation. We also obtain a degree-independent comparison for simultaneous approximation of finitely many $\varphi$-functions with a common denominator. Numerical pole distributions support the theoretical law.
\end{abstract}
\end{frontmatter}

\section{Introduction}

The convergence of rational approximations is closely tied to the distribution of their poles. In this paper, we study a converse question: to what extent does the approximation accuracy determine the pole distribution? Individual coefficients and exceptional poles may vary, but can the asymptotic distribution of almost all poles be uniquely determined by the convergence rate? We prove this implication for the exponential on an infinite sector. The optimal rate is known; the question concerns \emph{every} sequence of rational approximants attaining it.

We work with the following closed infinite sector in the Riemann sphere:
\begin{equation}\label{eq:sector-definition}
 \Sigma_\beta=\{re^{i\theta}:r\ge0,\ |\theta|\le\beta/2\}\cup\{\infty\},
 \qquad 0\le\beta<\pi.
\end{equation}
Here $\beta$ is the total opening and $\Chat=\C\cup\{\infty\}$. For $h>0$, define $f_h(z)=e^{-hz}$ at finite $z$ and $f_h(\infty)=0$. Write $\|g\|_E=\sup_{z\in E}|g(z)|$. The class $\RR_n(\Sigma_\beta)$ consists of rational functions $p/q$ with no pole in $\Sigma_\beta$ and, after cancellation, $\deg p\le\deg q\le n$. Poles are otherwise free outside the sector. Define
\begin{equation}\label{eq:error-definition}
 \eps_n(\beta;h)=\inf_{r\in\RR_n(\Sigma_\beta)}
                     \|f_h-r\|_{\Sigma_\beta}.
\end{equation}
Any rational function with finite error must be bounded on the positive real ray, hence has $\deg p\le\deg q$ after cancellation. Thus this is precisely the type-at-most-$(n,n)$ uniform approximation problem with finite error.

\subsection{The classical rate and the geometric question}
Gonchar and Rakhmanov (GR) determined the root limit in \eqref{eq:error-definition}, including the general infinite sector \cite[\S2.5, p.~329, Eq.~(27)]{GoncharRakhmanov1989}. We use a real-parameter form of their result, whose equivalence to GR's theta equation is proved in Proposition~\ref{prop:theta-ferrers}. Put
\[
 \eta=\frac{\pi-\beta}{2\pi},\qquad
 \mathsf P_\eta(1-2m_\beta)=0,\qquad \frac12<m_\beta<1,
\]
where $\mathsf P_\eta$ is the Legendre (Ferrers) function of the first kind. The root is unique, as shown in Lemma~\ref{lem:Legendre-root}. With the elliptic integral convention
\[
 \Kell(m)=\int_0^{\pi/2}\frac{dt}{\sqrt{1-m\sin^2t}},
\]
the classical result is
\begin{equation}\label{eq:intro-rho}
 \lim_{n\to\infty}\eps_n(\beta;h)^{1/n}
 =\rho_\beta:=\exp\!\left[-\pi\frac{\Kell(1-m_\beta)}{\Kell(m_\beta)}\right].
\end{equation}
Here $m$ is the elliptic \emph{parameter}, not the modulus. The limit is independent of $h$; for $\beta=0$ write $v=\rho_0=0.1076539192264845766\ldots$, the reciprocal Halphen constant. The parameter $m_\beta$ determines the classical rate and, after mass normalization, the shape and scale of the charged arc constructed in Section~\ref{sec:Scurve}.

We call a sequence $r_n\in\RR_n(\Sigma_\beta)$ \emph{root-optimal} if
\[
 \|f_h-r_n\|_{\Sigma_\beta}^{1/n}\longrightarrow\rho_\beta.
\]
This allows multiplicative error losses of size $e^{o(n)}$, without exact interpolation or optimality at each degree. Nevertheless, it determines the limiting pole measure.

\begin{theorem}[Pole rigidity of root-optimal sequences]\label{thm:root-optimal-rigidity}
For each $0\le\beta<\pi$ there is a probability measure $\nu_\beta$, supported on a compact analytic arc $\Gamma_\beta\subset\C\setminus\Sigma_\beta$, with the following property. If $h>0$ and $r_n\in\RR_n(\Sigma_\beta)$ is root-optimal, then
\begin{equation}\label{eq:all-root-optimal-poles}
 \frac1n\sum_{p\text{ a pole of }r_n}\delta_{hp/n}\weakto\nu_\beta.
\end{equation}
The sum counts poles with multiplicity \emph{after cancellation}; $\delta_a$ denotes unit mass at $a$, and convergence is weak convergence of measures on $\Chat$. No limiting mass lies at infinity or on the sector. The measure is the push-forward, under $z=-u^\gamma$ with $\gamma=2-\beta/\pi$, of the unit Green equilibrium charge constructed in Section~\ref{sec:Scurve}. In particular, the conclusion holds for every choice of best approximants.
\end{theorem}

Proposition~\ref{prop:degree} shows that root optimality forces the actual pole count to be $n-o(n)$, explaining the normalization by $n$. Weak convergence allows $o(n)$ exceptional poles whose locations need not follow the limiting arc; \ref{app:exceptional-poles} gives explicit constructions. Section~\ref{sec:Scurve} constructs the measure; Section~\ref{sec:minimax} proves convergence for every root-optimal sequence.

\subsection{What is inherited and what is extended}
The classical problem was motivated by numerical evolution equations \cite{CodyMeinardusVarga1969}; early sector estimates appear in \cite{SaffSchonhageVarga1976,SaffVargaNi1976}. GR obtained the optimal sector rate through a varying-weight equilibrium problem with equal side-normal derivatives, the \emph{S-property}. Our explicit realization supplies the positive charge, its scale, and the full variational inequality needed below; it does not replace GR's general approximation theorem.

Stahl--Schmelzer \cite[Theorems~9 and 10(i),(iii)]{StahlSchmelzer2009} stated the half-line pole law for \emph{close-to-best} sequences, defined by root optimality. Their article gives proof sketches and refers full proofs to a manuscript. Our principal extension is to $0<\beta<\pi$, within a proof that also includes the half-line. Schmelzer's recent paper \cite{Schmelzer2026} concerns half-line root rates for perturbed exponentials and $\varphi$-functions, using GR's varying-weight theorem.

Baratchart--Stahl--Yattselev (BSY) provide the closest general precedents. Their 2012 result treats irreducible critical points in the exterior Hardy space $\overline H_0^2$ through the reflected interpolation conditions supplied by criticality \cite[Theorem~7]{BaratchartStahlYattselev2012}. Their 2024 result gives pole rigidity and capacity error asymptotics for arbitrary root-optimal rational and meromorphic approximants in a fixed-target continuation class \cite[Theorems~2.1 and 2.5]{BaratchartStahlYattselev2024}.

The hypotheses of these results differ from ours. The point at infinity is an essential singularity \emph{on the boundary of the norming set}, so the exponential is not holomorphic in a neighborhood of that set. After degree scaling and a M\"obius change of coordinate (Section~\ref{sec:setup}), the target becomes $\exp[-n(1-w)/(1+w)]$ on a fixed compact lens. Lemma~\ref{lem:regularize} replaces this target near the lens by a holomorphic Cauchy transform $F_n$, but $F_n$ still varies exponentially with $n$ and carries a nonzero external field. Applying a fixed-target theorem separately for each $n$ does not justify a limit along the diagonal sequence of degree $n$.

The geometry alone is not the obstacle: the 2024 theorem assumes only a Jordan curve and allows lens corners when $\beta>0$. Open-arc Cauchy transforms may also have logarithmic endpoint branches. The decisive distinction for applying the 2024 fixed-target theorem is the exponential dependence on $n$; the 2012 theorem additionally assumes criticality and reflected interpolation, which root optimality alone does not provide.

Two obstacles remain for an arbitrary root-optimal sequence: poles may approach the compactified sector boundary, and approximation error gives a small Hermite residual rather than exact orthogonality. We first exclude positive pole mass on that boundary. For a hypothetical incorrect interior pole limit, we then adapt GR's strict-minimum, reflection, and noncancellation mechanisms \cite[\S3, Lemmas~8--9]{GoncharRakhmanov1989}. Balayage normalization and a rational phase correction give incompatible exponential bounds for \emph{the same Hermite integral}, while preventing complex cancellation in the lower bound. Spherical potentials avoid imposing logarithmic moments on the limiting pole measure. Thus root-optimal accuracy alone yields the pole law.

\subsection{Consequences and organization}
For positive-opening sectors, Section~\ref{sec:consequences} also treats the standard $\varphi$-functions used in exponential integrators. We prove that the best error for a fixed finite family sharing one denominator is comparable to the scalar best error, with constants independent of degree. Shared poles are familiar in exponential integrators and Carath\'eodory--Fej\'er approximation \cite{SchmelzerTrefethen2007,AlMohy2025}; the new result is this comparison on the whole infinite sector. A sublinear translation also selects root-optimal families whose poles all escape finite sets; this is not asserted for arbitrary or exact-best sequences.

Numerical candidates from adaptive Antoulas--Anderson (AAA) approximation followed by Lawson iteration \cite{NakatsukasaSeteTrefethen2018,NakatsukasaTrefethen2020}, fitted without prescribed pole geometry, test the predicted distribution in Section~\ref{sec:num}. Supplementary Material, Section S3, gives the additional matrix-function interpretation via the Crouzeix--Palencia estimate \cite{CrouzeixPalencia2017}.

The remainder of the paper is organized as follows. Section~\ref{sec:setup} introduces scaling and the Cauchy representation; Section~\ref{sec:Scurve} constructs the equilibrium charge and identifies the rate. Section~\ref{sec:minimax} proves rigidity in three stages, and Section~\ref{sec:consequences} develops consequences. Appendices~A--C supply auxiliary proofs; the independent supplement covers finite-degree bounds, wide-sector asymptotics, computation, and matrix functions.

\section{The approximation problem and its natural scaling}\label{sec:setup}

We first establish two finite-degree facts: best approximants exist and have positive error, and dilation removes the parameter $h$. We then rescale by the degree to obtain the varying target used in the equilibrium argument. This identifies the natural coordinate $hp/n$; the rigidity theorem proves that almost all poles occupy this scale. Additional finite-degree bounds appear in Supplementary Material, Section S1.

For later contour arguments, define the \emph{spherical order} of a reduced rational function $r=p/q$ by $\operatorname{ord}(r)=\max\{\deg p,\deg q\}$. It counts all poles on the Riemann sphere, including a possible pole at infinity.

\begin{proposition}[Existence and positivity at finite degree]
\label{prop:finite-best}
For each fixed $n\ge0$, $0\le\beta<\pi$, and $h>0$, the infimum $\eps_n(\beta;h)$ is attained and is strictly positive.
\end{proposition}

\begin{proof}
Choose a minimizing sequence $r_j=p_j/q_j$. Its sector norms are bounded, say by $M$, because the errors tend to $\eps_n$. Normalize the joint coefficient vectors of $(p_j,q_j)$ to norm one and pass to a polynomial limit $(p,q)$. At every finite $z\in\Sigma_\beta$,
\[
 |p(z)|\le M|q(z)|,\qquad |f_h(z)q(z)-p(z)|\le\eps_n|q(z)|.
\]
If $q\equiv0$, the first inequality forces $p$ to vanish on a real interval, contradicting the normalization. Otherwise, the first inequality makes every zero of $q$ in the sector removable after cancellation, and the second gives the required error at all finite points, including these removable zeros. Boundedness on the positive ray forces $\deg p\le\deg q$; letting $z\to\infty$ along that ray gives the error bound at infinity. Thus the minimum is attained. If it were zero, the identity theorem would make $e^{-hz}$ rational, which is impossible.
\end{proof}

\begin{proposition}[Scale invariance]\label{prop:scale}
For every $h>0$, $\eps_n(\beta;h)=\eps_n(\beta;1)$. Under the change of variable $\zeta=hz$, a pole $p$ becomes $hp$; hence $hp/n$ is the normalized pole coordinate in Theorem~\ref{thm:root-optimal-rigidity}.
\end{proposition}

\begin{proof}
The map $z\mapsto hz$ sends $\Sigma_\beta$ onto itself and preserves rational type. For $r(z)$ approximating $e^{-hz}$, set $s(\zeta)=r(\zeta/h)$. Then $s$ approximates $e^{-\zeta}$ with the same uniform error; the inverse substitution gives equality of the best errors.
\end{proof}

\subsection{The varying target and its Green equilibrium}\label{sec:external}

The root factor is already known from GR. To prove rigidity, we additionally need a Cauchy representation of the degree-scaled exponential and the equilibrium charge for its varying density. Two coordinate changes play different roles below: a M\"obius map preserves rational degree, whereas a fractional-power map is used only for the Green equilibrium.

Write $\mathfrak e_n(f;E)$ for the best type-at-most-$(n,n)$ uniform approximation error to $f$ on $E$. Since multiplication by $n$ maps $\Sigma_\beta$ onto itself without changing rational type, Proposition~\ref{prop:scale} gives
\begin{equation}\label{eq:varying-target}
 \eps_n(\beta;1)=\mathfrak e_n(e^{-nz};\Sigma_\beta).
\end{equation}
The degree and the exponential now involve the same parameter $n$. This identifies the natural pole coordinate for the analysis; only the rigidity theorem establishes the corresponding scale for arbitrary root-optimal sequences.

Fix $\beta/2<\vartheta<\pi/2$ and choose a contour $\mathcal F$ in the complement of the sector, positively oriented relative to $\Sigma_\beta$, whose two tails are eventually the rays $\arg t=\pm\vartheta$.  Then
\begin{equation}\label{eq:cauchy-inf}
 e^{-nz}=\frac{1}{2\pi i}\int_{\mathcal F}
 \frac{e^{-nt}}{t-z}\,dt,
 \qquad z\in\Sigma_\beta.
\end{equation}
Let $\mathcal F_R$ denote the contour obtained by truncating both tails at radius $R$, and let $f_{n,R}$ be its Cauchy integral. Then
\begin{equation}\label{eq:tail}
 \|e^{-nz}-f_{n,R}\|_{\Sigma_\beta}
 \le C_{\beta,\vartheta}\frac{e^{-nR\cos\vartheta}}{nR}.
\end{equation}
The density $e^{-nt}$ therefore induces the Gonchar--Rakhmanov external field
\begin{equation}\label{eq:field}
 \varphi(t)=\lim_{n\to\infty}\frac1{2n}
 \log\frac1{|e^{-nt}|}=\frac12\Rea t.
\end{equation}
With this convention, $|e^{-nt}|=e^{-2n\varphi(t)}$, and the equilibrium level gives the rational root factor $e^{-2\ell}$; this fixes the factor two in Theorem~\ref{thm:root}.

To compactify the norming set without changing rational type, we use a M\"obius coordinate. After fixing the truncated contour $\mathcal F_R$, choose $a_0>0$ with $-a_0\notin\mathcal F_R$ and put
\begin{equation}\label{eq:GR-Mobius}
 w=M(z)=\frac{a_0-z}{a_0+z},\qquad
 z=\Phi(w)=a_0\frac{1-w}{1+w}.
\end{equation}
Since positive dilation preserves $\Sigma_\beta$, its image $K_\beta=M(\Sigma_\beta)$ is independent of $a_0$: it is a compact lens (a segment when $\beta=0$). The points $w=1$ and $w=-1$ represent $z=0$ and $z=\infty$, respectively. Let $F_R=M(\mathcal F_R)$.  A direct change of variables gives
\begin{equation}\label{eq:GR-compact-Cauchy}
 f_{n,R}(\Phi(w))=c_{n,R}+\frac{1}{2\pi i}\int_{F_R}
       \frac{e^{-n\Phi(\tau)}}{\tau-w}\,d\tau,
\end{equation}
where
\[
 c_{n,R}=-\frac{1}{2\pi i}\int_{F_R}
          \frac{e^{-n\Phi(\tau)}}{1+\tau}\,d\tau.
\]
Indeed,
\[
 \frac{d\Phi(\tau)}{\Phi(\tau)-\Phi(w)}
 =\left(\frac1{\tau-w}-\frac1{1+\tau}\right)d\tau.
\]
The compact sets $K_\beta$ and $F_R$ are disjoint. On $F_R$, the Cauchy density has fixed analytic amplitude $1$ and the varying, nonvanishing factor $e^{-n\Phi}$; its external field is $\Rea\Phi/2$. The additive constant $c_{n,R}$ does not change rational type. This is the rectifiable-arc Cauchy setting of GR \cite[Theorem~1 and its arc formulation, p.~310]{GoncharRakhmanov1989}: the fixed function there is the Cauchy transform of the unit density on the arc, \emph{not} the constant function $1$. We will choose $R$ so that the truncation error is exponentially smaller than the optimal error.

We now switch from rational functions to the Green equilibrium. Unlike the M\"obius map, the following fractional-power map is applied only to measures and potentials, never to rational approximants. Let $D_\beta=\C\setminus\Sigma_\beta$ and $\gamma=2-\beta/\pi\in(1,2]$. Choose the branch
\begin{equation}\label{eq:power-map}
 u=(-z)^{1/\gamma},\qquad (-z)^{1/\gamma}>0\quad\text{for }z<0,
\end{equation}
which maps $D_\beta$ conformally onto the right half-plane.  Its Green kernel becomes
\begin{equation}\label{eq:green-RHP}
 g(u,v)=\log\left|\frac{u+\overline v}{u-v}\right|.
\end{equation}
For a positive measure $\nu$, write $G^\nu(u)=\int g(u,v)\,d\nu(v)$ for its Green potential. Since $z=-u^\gamma$, the external field $\varphi(z)=\Rea z/2$ becomes $-\Rea(u^\gamma)/2$. For probability measures supported on a candidate contour in the right half-plane, the weighted Green energy (Green interaction plus twice the external field) is
\begin{equation}\label{eq:green-energy}
 \mathcal J(\nu)
 =\iint g(u,v)\,d\nu(u)d\nu(v)
 -\int\Rea(u^\gamma)\,d\nu(u).
\end{equation}
Section~\ref{sec:Scurve} constructs a contour and positive unit charge satisfying the equilibrium and S-variational conditions for this field.

\section{An explicit realization of the classical sector equilibrium}\label{sec:Scurve}

We construct the GR sector equilibrium in four steps. A Ferrers zero fixes the angle of the two endpoints; a Stieltjes integral gives the candidate differential; a zero-level trajectory joins the endpoints; and its jump yields a positive charge whose total mass fixes the remaining scale. Proposition~\ref{prop:admissible} also verifies the global variational inequality. In particular, the trajectory is constructed before we identify it as the support of an equilibrium measure.

\subsection{The scalar parameter and the analytic differential}

Set $\eta=(\gamma-1)/2=(\pi-\beta)/(2\pi)\in(0,1/2]$.  For a conjugation-symmetric measure $\nu$ in the right half-plane, the derivative of its complex Green potential takes the form
\begin{equation}\label{eq:A}
 A(u)=-\int\frac{d\nu(v)}{u-v}+\int\frac{d\nu(v)}{u+v}.
\end{equation}
It is even, so write $A(u)=B(u^2)$. In this section $w=u^2$ denotes the squared half-plane coordinate, distinct from the M\"obius coordinate in \eqref{eq:GR-Mobius}. We first construct a candidate for this differential and then recover its positive charge. The angle $\theta$ below locates auxiliary branch points; it is not the sector half-opening. Take $\Lambda>0$ and $\pi/2<\theta<\pi$, and put
\[
 a_\pm=\sqrt\Lambda e^{\pm i\theta/2},\quad
 A_\pm=\Lambda e^{\pm i\theta},\quad c=\cos\theta,
 \quad R(w)^2=(w-A_+)(w-A_-).
\]
\begin{figure}[!hbp]
\centering
\begin{tikzpicture}[x=1cm,y=1cm,scale=.78]
\draw[->] (0,-3.05)--(0,3.05) node[above] {$\Ima u$};
\draw[->] (-.15,0)--(5.6,0) node[right] {$\Rea u$};
\coordinate (ap) at (2.65,1.35); \coordinate (am) at (2.65,-1.35); \coordinate (x0) at (1.45,0); \coordinate (iyb) at (0,2.25);
\draw[line width=1.05pt] (ap) .. controls (2.15,.95) and (1.65,.35) .. (x0) .. controls (1.65,-.35) and (2.15,-.95) .. (am);
\draw[dashed] (ap) .. controls (1.75,1.75) and (.65,2.12) .. (iyb);
\draw[dashed,->] (ap) .. controls (3.45,1.55) and (4.25,2.10) .. (5.15,2.72) node[right] {$\infty$};
\draw[dashed] (am) .. controls (1.75,-1.75) and (.65,-2.12) .. (0,-2.25);
\draw[dashed,->] (am) .. controls (3.45,-1.55) and (4.25,-2.10) .. (5.15,-2.72);
\fill (ap) circle (1.7pt) node[right] {$a_+$}; \fill (am) circle (1.7pt) node[right] {$a_-$}; \fill (x0) circle (1.5pt) node[above left] {$x_0$}; \fill (iyb) circle (1.5pt) node[left] {$iy_b$};
\node at (2.23,.12) {$\Gamma$};
\end{tikzpicture}
\caption{Topological zero-level geometry in the $u$-plane; angles and distances are not to scale. The solid arc $\Gamma$ is the charged support. The other upper trajectories reach $iy_b$ and infinity; reflection gives the lower half. Lemma~\ref{lem:trajectory} proves this outlet count. The initial circular cut in the $w=u^2$ plane is not the physical support.}\label{fig:S-curve-schematic}
\end{figure}
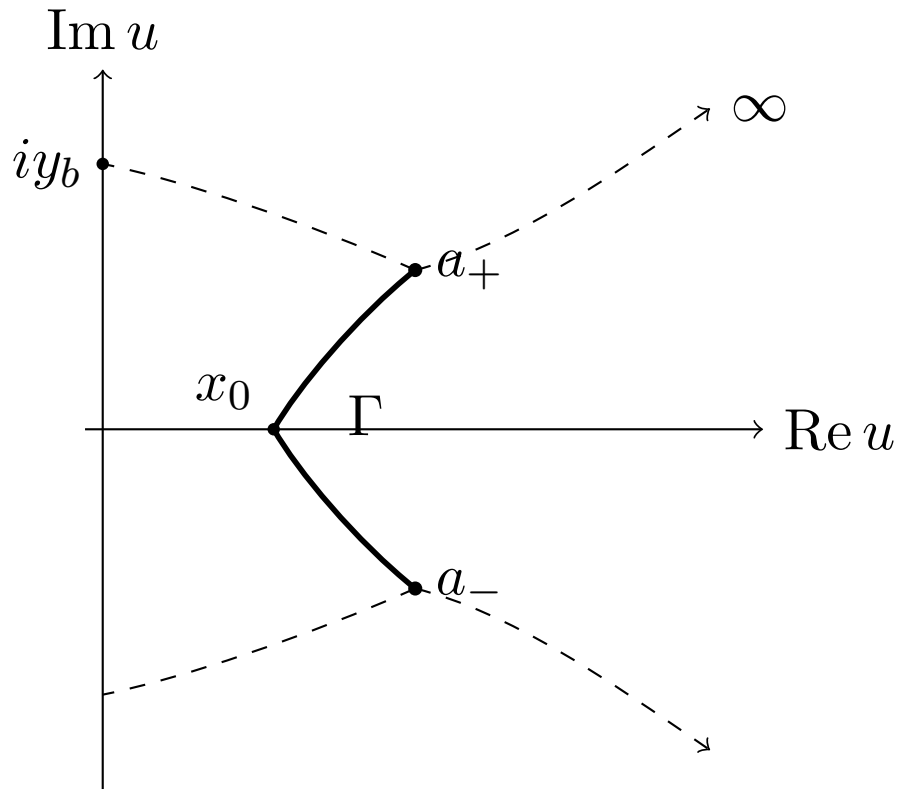
Choose the branch $R(w)\sim w$ at infinity. Initially the cut is the circular arc from $A_-$ to $A_+$ through the positive axis; after finding the equilibrium support, we move the cut to that arc. The same normalization gives $R(0)=-\Lambda$ for either cut. The required soft-edge jump and normalization conditions are
\begin{equation}\label{eq:Bjump}
 B_++B_-=\gamma w^\eta,
 \qquad B(w)=O(w^{-1})\quad(w\to\infty),
\end{equation}
with $B-\gamma w^\eta/2$ vanishing to square-root order at the two endpoints.  The constant-potential and S conditions require opposite boundary values of $A(u)-\gamma u^{\gamma-1}/2$. Thus $A_++A_-=\gamma u^{\gamma-1}$, which becomes the first relation in \eqref{eq:Bjump} when $w=u^2$. The square-root condition makes the endpoints soft edges rather than singular sources; without it, inverse-square-root homogeneous solutions would destroy uniqueness.

For an oriented circular temporary cut $\Gamma_0$ from $A_-$ to $A_+$, Mehler--Dirichlet's formula gives
\begin{equation}\label{eq:mehler}
 \int_{\Gamma_0}\frac{t^\eta}{R_+(t)}\,dt
       =-i\pi\Lambda^\eta\mathsf P_\eta(c).
\end{equation}
Indeed, set $t=\Lambda e^{is}$, $-\theta<s<\theta$, and use $R_+(t)=-\Lambda e^{is/2}\sqrt{2(\cos s-c)}$. Vanishing of this half-period gives the compatibility condition
\begin{equation}\label{eq:moment}
 \mathsf P_\eta(c)=0.
\end{equation}
The Stieltjes construction below separately checks the normalization at infinity: the period condition alone does not establish the required analytic solution. Thus \eqref{eq:moment} determines the endpoint angle, while the overall scale $\Lambda$ will be fixed by unit mass.

\begin{lemma}[Unique scalar parameter]\label{lem:Legendre-root}
For every $\eta\in(0,1/2]$, $\mathsf P_\eta$ has exactly one zero $c_\beta\in(-1,0)$.  With $m_\beta=(1-c_\beta)/2$, it is equivalently determined by
\[
 {}_2F_1(-\eta,1+\eta;1;m_\beta)=0,\qquad 1/2<m_\beta<1.
\]
\end{lemma}

\begin{proof}
Write $m=(1-c)/2$.  The Gauss hypergeometric expansion is $1-\sum_{k\ge1}b_km^k$, where
\[
 b_1=\eta(1+\eta)>0,\qquad
 b_{k+1}=b_k\frac{(k-\eta)(k+1+\eta)}{(k+1)^2}>0.
\]
It is strictly decreasing in $0<m<1$.  With $\Gamma(\cdot)$ denoting Euler's gamma function,
\[
 \mathsf P_\eta(0)=
 \frac{\sqrt\pi}{\Gamma((1-\eta)/2)\Gamma(1+\eta/2)}>0,
 \qquad b_k\sim\frac{\sin(\pi\eta)}{\pi k}.
\]
Hence $\sum_{k\ge1}b_k=\infty$, and since $b_km^k\uparrow b_k$ termwise as $m\uparrow1$, monotone convergence gives $1-\sum_{k\ge1}b_km^k\downarrow-\infty$.  Because $c=1-2m$, the endpoint $c=0$ corresponds to $m=1/2$ and $c\downarrow-1$ to $m\uparrow1$; strict decrease therefore gives exactly one zero with $c\in(-1,0)$, equivalently $m\in(1/2,1)$.  Our Ferrers and elliptic conventions follow \cite[Secs.~14.5, 14.12]{DLMF}.
\end{proof}

Hereafter $c=c_\beta$, unless another value is specified.  Define
\begin{equation}\label{eq:Stieltjes}
 S(w)=\int_0^\infty
 \frac{t^\eta\,dt}{\sqrt{t^2+2\Lambda ct+\Lambda^2}\,(t+w)},
 \qquad \kappa=\frac{\gamma\sin(\pi\eta)}{2\pi}>0,
\end{equation}
where the square root in the real integral is positive, and set
\begin{equation}\label{eq:BStieltjes}
 B(w)=\frac\gamma2w^\eta-\kappa R(w)S(w),\qquad
 F(u)=B(u^2)-\frac\gamma2u^{\gamma-1}.
\end{equation}
For this candidate define $A(u):=B(u^2)$. On the physical sheet, this gives
\begin{equation}\label{eq:F}
 F(u)=A(u)-\frac\gamma2u^{\gamma-1}.
\end{equation}
Proposition~\ref{prop:admissible} will verify that this candidate $A$ is indeed the Green derivative \eqref{eq:A} of the charge recovered from its jump.

The next lemma checks the Stieltjes normalization and establishes the first-order identity used to analyze the trajectories.

\begin{lemma}[Mellin normalization and the differential identity]
\label{lem:mellin-differential}
The function $B$ in \eqref{eq:BStieltjes} is analytic off the temporary cut, including across the negative real axis and at zero, and satisfies \eqref{eq:Bjump}.  Put $r(u)=R(u^2)$, so $F(u)=-\kappa r(u)S(u^2)$ on the physical sheet.  Then
\begin{equation}\label{eq:audited-ODE}
 \begin{split}
 &uF'(u)-(\gamma-1)F(u)=\frac{2C}{r(u)},\\
 &C=\frac{\gamma\eta}{2}\Lambda^{\eta+1}\mathsf P_{\eta-1}(c)
  =-\frac{\gamma(\eta+1)}2\Lambda^{\eta+1}\mathsf P_{\eta+1}(c)>0.
 \end{split}
\end{equation}
The identity also holds after replacing both $F$ and $r$ by their negatives.  In particular,
\begin{equation}\label{eq:F2}
 F(u)^2=\kappa^2(u^2-A_+)(u^2-A_-)S(u^2)^2
\end{equation}
has just the two simple zeros $a_\pm$ in the right half-plane.
\end{lemma}

The Stieltjes jump is first canceled by the fractional-power term. The combination $R(w)\{w\mathcal F'(w)-\eta\mathcal F(w)\}$, with $\mathcal F=B-\gamma w^\eta/2$, then has no jump and is constant by Liouville's theorem. Its value at zero is positive because $S(0)>0$. Details are in \ref{app:mellin}.

The mass and equilibrium level are determined by one half-period. To evaluate it with an unambiguous sign, use an auxiliary sheet specified by $\widehat r(0)=\Lambda$ and $\widehat F(u)=-\kappa\widehat r(u)S(u^2)$. Near zero, $\widehat r(u)=-R(u^2)$ and $\widehat F(u)=-F(u)$. The physical differential is still $F$ in \eqref{eq:F}.

\subsection{Periods, the connecting arc, and the positive charge}

\begin{lemma}[A branch-explicit half-period]\label{lem:half-period-audit}
Let $\widehat r(0)=\Lambda$ and continue $\widehat r(u)^2=(u^2-A_+)(u^2-A_-)$ from zero to $a_+$ in the first quadrant.  If $m=m_\beta=(1-c)/2$, then
\begin{equation}\label{eq:half-period-audit}
 I_a:=\int_0^{a_+}\frac{du}{\widehat r(u)}
       =\frac{\Kell(1-m)+i\Kell(m)}{2\sqrt\Lambda}.
\end{equation}
On the two-sheeted first-quadrant cover let $T(u)=\int_{a_+}^u\widehat F(s)\,ds$, with $\widehat F=-\kappa\widehat rS$ on this sheet.  Then
\begin{equation}\label{eq:T-audit}
 T(u)=\frac{u\widehat F(u)-2C\int_0^u ds/\widehat r(s)+2CI_a}{\gamma},
 \quad
 \Rea T(0)=\ell_0:=\frac{C\Kell(1-m)}{\gamma\sqrt\Lambda}>0.
\end{equation}
\end{lemma}

The substitution $u=\sqrt\Lambda\tan v$ turns the differential into an elliptic integral. Continuing first along a real segment and then vertically fixes the positive real and imaginary parts in \eqref{eq:half-period-audit}; integrating \eqref{eq:audited-ODE} gives \eqref{eq:T-audit}. The branch calculation is given in \ref{app:period}.

The half-period determines the level and the mass once the trajectory is known. The remaining question is geometric: which zero-level trajectory connects the two endpoints? The next lemma answers it by identifying the three outlets from the upper branch point.

\begin{lemma}[The connecting critical arc]\label{lem:trajectory}
There is a simple conjugation-symmetric critical arc $\Gamma$ in the right half-plane joining $a_-$ to $a_+$ and crossing the positive real axis exactly once.  Its defining equation on the cover is $\Rea T=0$.
\end{lemma}

\begin{proof}
The two-sheeted cover of the first quadrant branched at its single point $a_+$ is simply connected.  The differential $\widehat F(u)du$ is holomorphic there, with a zero of order two in the local square-root coordinate.  Its primitive $T$ is therefore single-valued; sheet exchange sends $T$ to $-T$.  Equivalently, the projected curves $\Rea T=0$ are critical trajectories of the quadratic differential $\widehat F(u)^2du^2$.  Locally $T(u)=c_*(u-a_+)^{3/2}(1+O(u-a_+))$, $c_*\ne0$, so three projected zero-level trajectories issue from $a_+$; see also Strebel \cite[Chs.~4--5]{Strebel1984}. We determine their endpoints by counting the available outlets.

Lemma~\ref{lem:mellin-differential} extends $B$ across the negative real axis. Hence $\widehat F=\pm(B(u^2)-\gamma u^{\gamma-1}/2)$ extends locally across each nonzero finite axis point and has a finite limit at zero. Its primitive $T$ is continuous up to these boundaries on each sheet.

On the positive real axis take $\widehat r(x)>0$.  Then $\partial_x\Rea T(x)=\widehat F(x)=-\kappa \widehat r(x)S(x^2)<0$, so $\Rea T(x)$ decreases from $\ell_0$ to $-\infty$ and has one transverse zero $x_0$.  On the positive imaginary axis, continue from zero with $\widehat r(iy)>0$.  The Stieltjes boundary values give
\[
 \Ima \widehat F(iy)=\frac\gamma2\sin(\pi\eta)y^{\gamma-1},\qquad
 \Rea T(iy)=\ell_0-\tfrac12\sin(\pi\eta)y^\gamma.
\]
Its derivative is $-\gamma\sin(\pi\eta)y^{\gamma-1}/2<0$, giving one transverse zero $iy_b$. At zero the two sheet limits have real parts $\pm\ell_0\ne0$, so zero is not an outlet.  Finally, on the sheet with $\widehat r(u)\sim u^2$ at infinity,
\[
 T(u)=-\tfrac12u^\gamma+T_\infty+O(u^{-1}).
\]
The angular derivative of its leading real part is nonzero near $\arg u=\pi/(2\gamma)$.  For every sufficiently large radius there is exactly one projected zero-level point, hence just one outlet at infinity.  The other sheet gives the same projected level set.

Three projected trajectories issue from $a_+$.  Along each oriented lifted regular trajectory $T$ is purely imaginary and $dT=\widehat F(u)du$ never vanishes, so $\Ima T$ is strictly monotone.  A branch therefore cannot return to $a_+$ or close up.  Nor can it accumulate in a compact regular part of the cover: if a maximal branch had such an accumulation point, monotonicity would give a finite limit of $\Ima T$, while the local uniqueness of the regular level curve through that point would continue the branch past the alleged endpoint.  Thus a branch must end at a critical point, a boundary outlet, or infinity.  There is no second critical point in the first-quadrant cover.  A regular boundary zero admits only one incident branch, and the large-radius description admits only one branch at infinity.  Hence the three branches use the three distinct outlets $x_0$, $iy_b$, and infinity.

The three outlets are now occupied, and zero is excluded. Any remaining compact component therefore lies in the interior of the cover. It would contain a closed regular level curve, which bounds a region on the simply connected cover; the maximum principle would force $\Rea T$ to vanish there, an impossibility.  In particular the branch joining $a_+$ to $x_0$ is simple.  Reflection across the real axis supplies its lower half and proves the assertion.
\end{proof}

Having identified the connecting arc, move the physical cut of $R(u^2)$ to $\Gamma$ and its negative reflection. Orient $\Gamma$ from $a_-$ to $a_+$ and define
\begin{equation}\label{eq:nu-density}
 d\nu(u)=\frac{F_+(u)}{\pi i}\,du.
\end{equation}
Here the plus side is the left side of the oriented arc, not automatically the sheet used in \eqref{eq:T-audit}.

The trajectory is not yet an equilibrium support: its jump must define a positive measure with the correct Green potential and a global variational inequality. The next proposition verifies these properties, using the half-period to determine the mass and then fixing $\Lambda$ by unit normalization.

\begin{proposition}[Positive charge and admissible S-contour]\label{prop:admissible}
The measure \eqref{eq:nu-density} is positive, has square-root vanishing density at the endpoints, and has mass
\begin{equation}\label{eq:audited-mass}
 M=\frac{2C\Kell(m)}{\pi\gamma\sqrt\Lambda}.
\end{equation}
There is a unique $\Lambda>0$ for which $M=1$.  For that scale its Green potential in the field $-\Rea u^\gamma/2$ has level $\ell=\ell_0$ on $\Gamma$, satisfies the S-property, and the arc has two admissible tails on which the total potential is strictly larger than $\ell$.
\end{proposition}

\begin{proof}
We verify positivity and mass first, reconstruct the Green derivative from the jump, and finally construct uncharged tails satisfying the strict inequality.

\emph{Positivity and mass.} At the real crossing, $F_+(x_0)=-\widehat F(x_0)>0$ by the auxiliary-sheet calculation in Lemma~\ref{lem:half-period-audit}.  The tangent of the upward-oriented arc is positive imaginary there.  Thus \eqref{eq:nu-density} is positive at the crossing.  The trajectory condition makes it real along the arc; it cannot change sign before an endpoint, since $F$ has no further zeros.  Its endpoint order follows from \eqref{eq:F2}.

On the auxiliary sheet, $\Ima T(x_0)=C\Kell(m)/(\gamma\sqrt\Lambda)$. Consequently the mass on the upper half of $\Gamma$ is $\Ima T(x_0)/\pi$; reflection doubles it and gives \eqref{eq:audited-mass}. Since $C$ is a positive constant times $\Lambda^{\eta+1}$, the mass is a positive constant times $\Lambda^{\gamma/2}$, proving uniqueness of its unit normalization.

\emph{Reconstruction of the Green derivative.} The density was defined from a jump of $F$, so one must still check that no extra analytic term is present.  Apply Plemelj's formula to the Cauchy transform of \eqref{eq:nu-density} on $\Gamma$ and its negative reflection.  Its jump is exactly the jump prescribed in \eqref{eq:A}.  Subtracting this Cauchy transform from the explicit $A$ therefore removes both cuts; the square-root endpoint behavior makes the endpoint singularities removable, and the normalization at infinity makes the difference vanish there.  Liouville's theorem then forces the difference to be zero.  Thus the differential constructed above is the actual complex Green derivative of the positive charge, not merely a stationary formal solution.

Write $W=G^\nu-\Rea u^\gamma/2$.  On the physical sheet, \eqref{eq:A} and \eqref{eq:F} give $dW=\Rea\{F(u)\,du\}$.  In the first-quadrant component attached to zero, $F=-\widehat F$ and $dT=\widehat F\,du$.  Since $W(0)=0$ and $\Rea T(0)=\ell_0$, we obtain $W=\ell_0-\Rea T$ there.  Hence $W(x_0)=\ell_0$, and the zero-level trajectory $\Rea T=0$ is precisely the level $W=\ell_0$.  Thus $W=\ell_0$ on the entire charged arc.  The two side derivatives are opposite analytic branches, which is exactly the S-property.

\emph{Global inequality and admissible tails.} It remains to ensure that the charged arc is not only stationary but globally extremal on one admissible Cauchy contour.  For this, one needs the full trajectory graph, not merely an unboundedness assertion from the maximum principle.  In the upper quadrant the other two zero-level branches join $a_+$ to $iy_b$ and to infinity.  The tripod described in Lemma~\ref{lem:trajectory} bounds a component incident to $a_+$ whose outer boundary includes $iy$ with $y>y_b$ and the angles $\pi/(2\gamma)<\arg u<\pi/2$ at infinity.  It contains no additional zero level.  On its imaginary boundary
\[
 W(iy)-\ell_0=\tfrac12\sin(\pi\eta)y^\gamma-\ell_0>0,
\]
so this component is the positive one.  Moreover $G^\nu(u)=O(|u|^{-1})$, so outside a sufficiently large circle any fixed closed angular subsector of these angles lies in it.  A simple tail can therefore be drawn from $a_+$ through this same positive component and made eventually a fixed ray; reflect it for the lower tail.  The endpoint $3/2$-power local structure permits a rectifiable, piecewise analytic start. Under $z=-u^\gamma$ the eventual rays satisfy $\beta/2<|\arg z|<\pi/2$, as required by the Cauchy representation.

The resulting full contour has $W=\ell$ on $\supp\nu$ and $W>\ell$ on its uncharged parts. To see why this gives equilibrium, let $\widetilde\nu$ be any competing probability measure of finite Green energy on a finite truncation. Using \eqref{eq:green-energy} and $W=G^\nu-\Rea(u^\gamma)/2$,
\[
 \mathcal J(\widetilde\nu)-\mathcal J(\nu)
 =\iint g(u,v)\,d(\widetilde\nu-\nu)(u)d(\widetilde\nu-\nu)(v)
   +2\int (W-\ell)\,d\widetilde\nu\ge0.
\]
The first term is the nonnegative Green energy of the zero-mass difference, and the second is nonnegative because $W\ge\ell$ on the chosen contour.  This proves equilibrium on every finite truncation, hence the global S-capacitor inequalities on the chosen contour.
\end{proof}

\begin{proposition}[Recovery of the original capacitor charge]\label{prop:balayage}
The Green equilibrium measure determines a positive unit charge on the fixed sector plate, and hence the logarithmic S-capacitor of the original problem. This identifies the construction with GR's two-plate formulation; the rigidity proof uses its Green form directly.
\end{proposition}

\begin{proof}
For $v$ in the right half-plane, harmonic measure on the imaginary axis has the explicit Poisson density
\[
 d\omega_v(iy)=\frac{\Rea v}{\pi|iy-v|^2}\,dy.
\]
Average against $d\nu(v)$ and push forward by $z=-(iy)^\gamma$ to obtain $\nu_E$.  The measure is positive and has mass one.  Its density is $O(y^{-2})$ at infinity; after the power map this becomes $O(r^{-1-1/\gamma})$, and near zero it is $O(r^{1/\gamma-1})$. These bounds give finite logarithmic moments and energy and no atom at infinity.

Let $\nu_F$ be the push-forward of $\nu$ under $z=-u^\gamma$, and write $V^\alpha(z)=\int\log|z-t|^{-1}\,d\alpha(t)$.  Conformal invariance of harmonic measure and the Green kernel identifies $\nu_E$ as the harmonic balayage of $\nu_F$ onto the original sector.  Apply the point-charge Green representation in $D_\beta$ and then average against $\nu_F$.  The difference has zero boundary values and vanishes at infinity; the moment bounds above fix the logarithmic normalization.  Thus
\[
 V^{\nu_F-\nu_E}=G_{D_\beta}^{\nu_F}\quad\hbox{in }D_\beta,
 \qquad V^{\nu_F-\nu_E}=0\quad\hbox{on }\Sigma_\beta.
\]
Thus the Green equilibrium level and S-property correspond to the two-plate logarithmic formulation used by GR. The fractional-power map changes only measures and potentials; it is never applied to rational approximants, which remain in the original variable or a M\"obius coordinate.
\end{proof}

\subsection{The classical root factor and its normalization}\label{sec:root}

The preceding construction gives a positive unit charge and its equilibrium level, but it has not yet identified that charge with poles of arbitrary approximants. We first recover the classical root factor from the level and reconcile its Ferrers and GR theta descriptions; pole rigidity is proved in Section~\ref{sec:minimax}.

In addition to $\Kell$, define $\Eell(m)=\int_0^{\pi/2}(1-m\sin^2t)^{1/2}\,dt$, again using the elliptic parameter.
Lemma~\ref{lem:half-period-audit} and Proposition~\ref{prop:admissible} give, for a measure of mass $M$, $M=2C\Kell(m)/(\pi\gamma\sqrt\Lambda)$, where $m=(1-c)/2$, and equilibrium level $\ell=C\Kell(1-m)/(\gamma\sqrt\Lambda)$.
After imposing $M=1$ the scale cancels and
\begin{equation}\label{eq:ell}
 \boxed{\displaystyle
 \ell=\frac\pi2\frac{\Kell(1-m_\beta)}{\Kell(m_\beta)}.}
\end{equation}
The same branch-explicit half-period determines the support scale:
\begin{equation}\label{eq:lambda-scale}
 \Lambda^{\gamma/2}
 =\frac{\pi}{-(\eta+1)\mathsf P_{\eta+1}(c_\beta)\Kell(m_\beta)}.
\end{equation}
Thus the endpoint scale enters the predicted pole locations but not the root factor.

\begin{theorem}[GR sector root factor, in Ferrers form]\label{thm:root}
For every $0\le\beta<\pi$ and $h>0$,
\begin{equation}\label{eq:rho-main}
 \boxed{\displaystyle
 \lim_{n\to\infty}\eps_n(\beta;h)^{1/n}
 =\rho_\beta
 =\exp\!\left[-\pi\frac{\Kell(1-m_\beta)}{\Kell(m_\beta)}\right]
 =e^{-2\ell},}
\end{equation}
where $\ell$ is the level of the unit equilibrium charge in Proposition~\ref{prop:admissible}, and $m_\beta\in(1/2,1)$ is determined by
\[
 \mathsf P_{(\pi-\beta)/(2\pi)}(1-2m_\beta)=0.
\]
\end{theorem}

\begin{proof}
GR's sector result \cite[\S2.5, p.~329, Eq.~(27)]{GoncharRakhmanov1989} gives the root limit and its theta equation. Proposition~\ref{prop:theta-ferrers} below identifies that equation with the Ferrers equation, whose unique root is supplied by Lemma~\ref{lem:Legendre-root}. The elliptic nome relation gives \eqref{eq:rho-main}, and Proposition~\ref{prop:scale} removes $h$. Finally, \eqref{eq:ell} identifies the same factor as $e^{-2\ell}$ for the positive equilibrium charge just constructed. Thus GR supplies the root factor, while the explicit charge and its full variational conditions supply the rigidity input.
\end{proof}

\subsubsection*{The theta equation and the half-line calibration}
GR uses a half-opening $\theta_{\rm GR}=\beta/2$ and the exponent parameter $\alpha_{\rm GR}=1-\beta/(2\pi)=\eta+1/2$. More subtly, its theta nome is imaginary. The following identity fixes the correspondence, including its normalization.

\begin{proposition}[Theta--Ferrers correspondence]\label{prop:theta-ferrers}
Use the standard theta functions with nome $q$, so that
\[
 \theta_3(v,q)=\sum_{k\in\mathbb Z}q^{k^2}e^{2ikv},\qquad
 \theta_4(v,q)=\sum_{k\in\mathbb Z}(-1)^kq^{k^2}e^{2ikv},
\]
and let $\theta_2(v,q)=2\sum_{k\ge0}q^{(k+1/2)^2}\cos((2k+1)v)$ for positive $q$. For $q\in(0,1)$ set
\[
 m(q)=\frac{\theta_2(0,q)^4}{\theta_3(0,q)^4},\qquad
 q=\exp\!\left[-\pi\frac{\Kell(1-m(q))}{\Kell(m(q))}\right].
\]
Then, for $0<\eta\le1/2$,
\begin{equation}\label{eq:theta-ferrers}
 \int_0^{1/2}\left(
 \frac{\theta_4(\pi t,i\sqrt q)}{\theta_3(\pi t,i\sqrt q)}
 \right)^{2\eta+1}dt
 =\frac{\pi}{4\Kell(m(q))}\mathsf P_\eta(1-2m(q)).
\end{equation}
The power is defined by the principal logarithm: the theta quotient stays on the open right half of the unit circle. With GR's nome $\mathfrak h=i\sqrt q$, its root factor is $-\mathfrak h^2=q$. Hence its equation~(27) is equivalent to the Ferrers zero equation in Theorem~\ref{thm:root}.
\end{proposition}

\ref{app:theta} proves the identity by separating even and odd theta indices and applying the Mehler--Dirichlet formula. It holds away from the root as well. Thus the two descriptions specify exactly the same constant, not merely constants with agreeing numerical values.

At $\beta=0$, $\eta=1/2$ and
$\mathsf P_{1/2}(1-2m)=2\{2\Eell(m)-\Kell(m)\}/\pi$.
The equation therefore reduces to $\Kell(m)=2\Eell(m)$, giving $v=\rho_0$.
\ref{app:theta} also checks the endpoint scale against GR's half-line construction.
The wide-sector expansion and stable evaluation of $m_\beta$ are independent of the rigidity proof and are given in Supplementary Material, Section S2.

Define $\nu_\beta$ as the push-forward under $z=-u^\gamma$ of the positive unit measure \eqref{eq:nu-density}, with $\Lambda$ fixed by \eqref{eq:lambda-scale}. This measure is defined entirely by the S-capacitor construction, independently of the denominator of any approximant. Section~\ref{sec:minimax} proves that it is the limiting pole law of every reduced root-optimal sequence.

\section{Pole rigidity: from root-optimality to the equilibrium law}
\label{sec:minimax}

Theorem~\ref{thm:root-optimal-rigidity} concerns arbitrary root-optimal approximants, whose denominators need not arise from interpolation or orthogonality. We therefore work with a general subsequential pole limit, without assuming interpolation equations for the denominators. The proof has three stages, summarized by the roles they play in the final contradiction.
\begin{center}
\begin{tabular}{@{}p{.25\linewidth}p{.69\linewidth}@{}}
\toprule
Stage & Purpose and principal result\\
\midrule
Locate the pole mass & Degree saturation and pole stripping rule out loss of limiting pole mass at the compactified boundary (Proposition~\ref{prop:degree}, Theorem~\ref{thm:boundary}).\\[2pt]
Test a wrong limit & Green reflection produces a strict minimum with exponent $m_*<2\ell$ (Theorem~\ref{thm:green-GL9}, Lemma~\ref{lem:green-gap}).\\[2pt]
Compare one integral & A normalized zero table and a phase correction give root size $e^{-m_*}$; the approximation error bounds the same Hermite integral by $e^{-2\ell}$ (Theorem~\ref{thm:adaptive-sync}, Proposition~\ref{prop:sync-residual-upper}).\\
\bottomrule
\end{tabular}
\end{center}
We use the equilibrium constructed in Section~\ref{sec:Scurve}: its charge is a positive probability measure on a simple nonseparating analytic arc, its Green potential is continuous, and the full S-variational conditions hold. Its level $\ell$ satisfies $\rho_\beta=e^{-2\ell}$. In compact coordinates we denote the charge by $\lambda_F$ and write $\ell_\beta=\ell$. The Cauchy contour is the S-contour of Proposition~\ref{prop:admissible}, truncated beyond the charged arc, with varying density from Section~\ref{sec:external}.

\subsection{Coordinates and the degree budget}

Three coordinate operations appear in this section, and they have different purposes. The substitution $\zeta=hz/n$ sets the pole scale; M\"obius maps make the norming set compact without changing rational order. The power map from Section~\ref{sec:external} and the conformal maps used below act only on Green potentials and measures, never on the rational functions in the Hermite identity.

Let $r_n$ be any reduced root-optimal sequence; actual best approximants are a special case.  The natural scaled variable is $\zeta=hz/n$; then $r_n((n/h)\zeta)$ approximates $e^{-n\zeta}$ on the same sector with the same error.  Only after this scaling do we apply
\begin{equation}\label{eq:Mobius}
 w=\frac{1-\zeta}{1+\zeta},
 \qquad
 \zeta=\frac{1-w}{1+w}.
\end{equation}
The sector becomes the fixed lens $K_\beta=\{(1-\zeta)/(1+\zeta):\zeta\in\Sigma_\beta\}$, and the target and transformed approximant are
\begin{equation}\label{eq:Gn}
 G_n(w)=\exp\!\left[-n\frac{1-w}{1+w}\right],
 \qquad
 \mathcal R_n(w)=r_n\!\left(\frac{n}{h}\frac{1-w}{1+w}\right).
\end{equation}
Thus poles escaping to infinity on the scale $hp/n$ approach $w=-1$, while poles collapsing to zero approach $w=1$. Poles approaching either boundary ray of the scaled sector likewise approach $\partial K_\beta$. These are the boundary losses that the next theorem excludes.

Let $\Omega_\beta=\Chat\setminus K_\beta$, let $g_\beta(w,a)$ be its Green function, and define the Green height
\begin{equation}\label{eq:height}
 a_\beta(a)=g_\beta(a,\infty).
\end{equation}
The level curves $\Gamma_s=\{a_\beta=s\}$, $s>0$, are analytic Jordan curves.

\paragraph{Why almost all degrees must be used}

Before studying their locations, we must show that a root-optimal sequence actually uses asymptotically all available denominator degrees. Otherwise it would achieve the $n$-degree root rate with strictly fewer degrees.

\begin{proposition}[Degree usage]\label{prop:degree}
If $r_n=p_n/q_n\in\RR_n(\Sigma_\beta)$ is reduced and root-optimal, then its actual denominator degree $m_n=\deg q_n$ satisfies $m_n/n\to1$.
\end{proposition}

\begin{proof}
Finite uniform error on the unbounded sector implies $\deg p_n\le\deg q_n=m_n$, so $r_n$ belongs to the type-$(m_n,m_n)$ class.  Hence
\[
 \|e^{-hz}-r_n\|_{\Sigma_\beta}\ge\eps_{m_n}(\beta;h).
\]
If $m_n$ remains bounded along a subsequence, positivity of the fixed-degree best error from Proposition~\ref{prop:finite-best} makes its $n$th-root lower bound tend to one, a contradiction.  Thus $m_n\to\infty$.  If $m_n/n\to\alpha<1$ along a subsequence, Theorem~\ref{thm:root} gives
\[
 \rho_\beta\ge\rho_\beta^\alpha>\rho_\beta,
\]
a contradiction.
\end{proof}

\subsection{Boundary localization by pole stripping}

The next step rules out a positive fraction of poles approaching the compactified boundary. We find a Green level curve close to that boundary on which the approximation error remains exponentially small, then use Cauchy projection to remove poles on the bounded side of the curve. If a positive fraction were removed, the resulting lower-degree approximant would be too accurate for its degree, contradicting Theorem~\ref{thm:root}.

\begin{lemma}[Direct logarithmic-radius contour selection]\label{lem:cap-select}
Let $K\subset\C$ be a regular continuum with simply connected complement, $\Gamma_s=\{g_K(\cdot,\infty)=s\}$, and $0<s_1<s_2$.  For any multiset $p_1,\ldots,p_N$ in $\Chat\setminus K$, including possible points at infinity, with $N\le n$ and $n\ge1$, there exists $s\in[s_1,s_2]$ such that
\begin{equation}\label{eq:log-radius-selection}
 \frac1n\max_{z\in\Gamma_s}\sum_{j=1}^N g_K(z,p_j)
 \le\frac{s_2+\pi^2/2}{\log(s_2/s_1)}.
\end{equation}
In particular the selected curve contains no pole, and the bound can be made arbitrarily small by choosing a fixed $s_1$ sufficiently small relative to $s_2$, independently of $n$ and of all pole locations.
\end{lemma}

The estimate comes from averaging the explicit disk Green kernel with respect to the logarithmic radius $ds/s$. Uniformity in pole locations permits contour selection even for clusters, without residue estimates. The proof is in \ref{app:contour}.

\begin{lemma}[Cauchy pole stripping]\label{lem:pole-strip}
Let $\Gamma$ be a positively oriented Jordan curve and let $R$ be rational with no pole on $\Gamma$.  For $z$ in the bounded component of $\Chat\setminus\Gamma$ set
\[
 \Pi_\Gamma R(z)=\frac1{2\pi i}\int_\Gamma\frac{R(t)}{t-z}\,dt.
\]
Then $\Pi_\Gamma R$ deletes the principal parts at all finite poles in the bounded component of $\Gamma$.  It is rational, pole-free there, and its spherical order drops by at least the deleted multiplicity.
\end{lemma}

\begin{proof}
Apply the residue theorem to the integrand $R(t)/(t-z)$.  A pole of $R$ inside $\Gamma$ contributes exactly the negative of its principal part evaluated at $z$, whereas poles outside $\Gamma$ do not contribute.  The polynomial part, or equivalently a possible pole at infinity, is unchanged.  The spherical-order claim follows from this partial-fraction decomposition.
\end{proof}

\begin{lemma}[Holomorphic regularization of the varying target]\label{lem:regularize}
Fix $\kappa_{\mathrm{tail}}>2\ell$.  There are $T_0>0$ and functions $F_n$ holomorphic in a fixed neighborhood of the closed region bounded by $\Gamma_{T_0}$ such that
\begin{equation}\label{eq:Fn-close}
 \|F_n-G_n\|_{K_\beta}\le C e^{-\kappa_{\mathrm{tail}} n}
\end{equation}
and
\begin{equation}\label{eq:Fn-growth}
 \max_{\Gamma_{T_0}}|F_n|\le e^{C_0n}
\end{equation}
for constants $C,C_0$ independent of $n$.
\end{lemma}

\begin{proof}
Use the Cauchy contour in Section~\ref{sec:external} and truncate its two tails at a fixed radius $R$ so large that $R\cos\vartheta>\kappa_{\mathrm{tail}}$.  The resulting finite Cauchy transform $f_{n,R}$ satisfies the uniform tail estimate \eqref{eq:tail}. It is holomorphic off the fixed truncated contour and at infinity, where its value is zero. Define $F_n$ by transporting $f_{n,R}$ through the scale-first M\"obius map \eqref{eq:Mobius}. The transformed contour is a spherical compact set disjoint from $K_\beta$; hence some closed Green collar through $\Gamma_{T_0}$ avoids it. In this preliminary lens chart the contour may contain the point at infinity, but $F_n$ only needs to be holomorphic near the collar. Section~\ref{subsec:adaptive-sync} later chooses a finite Cauchy chart for the Hermite identity. The tail bound gives \eqref{eq:Fn-close}.  On the fixed truncated contour the density $e^{-nt}$ is bounded by $e^{C_0n}$, while the Cauchy kernel is uniformly bounded for points on $\Gamma_{T_0}$; hence \eqref{eq:Fn-growth}.
\end{proof}

\begin{theorem}[Boundary tightness of root-optimal poles]\label{thm:boundary}
Let $r_n\in\RR_n(\Sigma_\beta)$ be root-optimal.
Let $w_{j,n}$ be the poles, counted with multiplicity, after the scale-first transformation $\zeta=hz/n$ and the M\"obius map \eqref{eq:Mobius}.  Put
\[
 \mu_n=\frac1n\sum_j\delta_{w_{j,n}}.
\]
By Proposition~\ref{prop:degree}, $\mu_n(\Omega_\beta)\to1$.  Then
\begin{equation}\label{eq:boundary-tight}
 \boxed{\displaystyle
 \lim_{s\downarrow0}\ \limsup_{n\to\infty}
 \mu_n\{p:a_\beta(p)<s\}=0.}
\end{equation}
In particular no subsequential weak limit of the normalized pole measures places positive mass on $\partial K_\beta$.
\end{theorem}

\begin{proof}
We regularize the varying target at exponentially negligible cost, remove the residual pole singularities from a subharmonic estimate, select a Green level curve with small total pole potential, and strip the poles inside that curve. If there were at least $\delta n$ such poles, the reduced degree would contradict the optimal root rate.

Fix $0<\delta<1$. We shall find $s_\delta>0$ such that, for all sufficiently large $n$, fewer than $\delta n$ poles have Green height below $s_\delta$. Choose a number $\xi$ satisfying
\begin{equation}\label{eq:xi-choice}
 0<\xi<2\ell\delta,
\end{equation}
where $\rho_\beta=e^{-2\ell}$.  We work in the fixed lens variable.

Choose $\kappa_{\mathrm{tail}}>2\ell$ and let $F_n$ and $\Gamma_{T_0}$ be supplied by Lemma~\ref{lem:regularize}.  By root optimality and \eqref{eq:Fn-close},
\begin{equation}\label{eq:nearbest-F}
 \|F_n-\mathcal R_n\|_{K_\beta}=e^{-2\ell n+o(n)}.
\end{equation}
Moreover $\|\mathcal R_n\|_{K_\beta}\le1+o(1)$.

Let $p_{j,n}$ denote all poles of $\mathcal R_n$ on the sphere, repeated according to multiplicity, and set
\begin{equation}\label{eq:un}
 u_n(z)=\log|F_n(z)-\mathcal R_n(z)|-
 \sum_jg_\beta(z,p_{j,n}).
\end{equation}
The reason for subtracting the Green sum is local: at each pole of $\mathcal R_n$, its logarithmic singularity cancels the one in $\log|F_n-\mathcal R_n|$. Thus $u_n$ extends subharmonically through these poles in $0<a_\beta<T_0$. We can then apply the two-constants theorem without knowing their residues.  On the inner boundary,
\begin{equation}\label{eq:un-inner}
 u_n\le-2\ell n+o(n).
\end{equation}
On $\Gamma_{T_0}$, the standard rational Bernstein--Walsh inequality (see, e.g., \cite{Walsh1969}) gives
\[
 \log|\mathcal R_n(z)|-\sum_jg_\beta(z,p_{j,n})
 \le\log\|\mathcal R_n\|_{K_\beta}=O(1),
\]
while \eqref{eq:Fn-growth} and nonnegativity of the Green sum give $u_n\le C_1n$.  The harmonic measure of the inner boundary in this Green annulus is $1-a_\beta(z)/T_0$. The two-constants theorem therefore gives, uniformly for $z\in\Gamma_s$,
\begin{equation}\label{eq:two-constants}
 u_n(z)
 \le -2\ell n+\frac{s}{T_0}(C_1+2\ell)n+o(n).
\end{equation}
Choose $s_2>0$ so small that the added term in \eqref{eq:two-constants} is at most $\xi n/4$ for $0<s\le s_2$ and all large $n$.

Choose $L>4(s_2+\pi^2/2)/\xi$ and put $s_1=s_2e^{-L}$. Lemma~\ref{lem:cap-select} gives, for every $n$, a curve $\Gamma_{t_n}$ with $t_n\in[s_1,s_2]$ and
\[
 \max_{\Gamma_{t_n}}\sum_jg_\beta(z,p_{j,n})<\xi n/4.
\]
The heights $s_1,s_2$ are fixed independently of $n$, although the selected height $t_n$ may vary with $n$. Combining \eqref{eq:un} and \eqref{eq:two-constants} gives
\begin{equation}\label{eq:error-contour}
 \max_{\Gamma_{t_n}}|F_n-\mathcal R_n|
 \le \exp\{-2\ell n+\tfrac12\xi n+o(n)\}.
\end{equation}

Apply Lemma~\ref{lem:pole-strip} to $\mathcal R_n$ on $\Gamma_{t_n}$, and call the result $\widetilde R_n$. It deletes all poles on the bounded side of the curve, without requiring any residue estimate. For $z$ inside $\Gamma_{t_n}$,
\[
 \widetilde R_n(z)=\frac1{2\pi i}
 \int_{\Gamma_{t_n}}\frac{\mathcal R_n(t)}{t-z}\,dt.
\]
Since $F_n$ is holomorphic inside the contour,
\begin{equation}\label{eq:cauchy-strip}
 F_n(z)-\widetilde R_n(z)=\frac1{2\pi i}
 \int_{\Gamma_{t_n}}\frac{F_n(t)-\mathcal R_n(t)}{t-z}\,dt,
 \qquad z\in K_\beta.
\end{equation}
The family $\{\Gamma_s:s\in[s_1,s_2]\}$ has uniformly bounded length and stays a positive distance from $K_\beta$.  Hence \eqref{eq:error-contour} and \eqref{eq:cauchy-strip} imply
\begin{equation}\label{eq:stripped-error}
 \|F_n-\widetilde R_n\|_{K_\beta}
 \le \exp\{-2\ell n+\xi n+o(n)\}.
\end{equation}

Assume, to the contrary, that infinitely many $n$ have at least $\delta n$ poles with $a_\beta(p)<s_1$.  All of them lie in the bounded side of $\Gamma_{t_n}$.  Lemma~\ref{lem:pole-strip}, with poles counted spherically and with multiplicity, therefore gives
\[
 \operatorname{ord}(\widetilde R_n)\le n-\delta n.
\]
Thus, after replacing the right-hand side by its integer part, $\widetilde R_n$ has order at most
\[
 M_n=\lfloor(1-\delta)n\rfloor.
\]
Returning through the M\"obius map and adding the exponentially smaller regularization error \eqref{eq:Fn-close}, we obtain an admissible rational approximant of order at most $M_n$ to $e^{-nz}$ on $\Sigma_\beta$ satisfying
\[
 \eps_{M_n}(\beta;n)
 \le\exp\{-2\ell n+\xi n+o(n)\}.
\]
By scale invariance $\eps_{M_n}(\beta;n)=\eps_{M_n}(\beta;1)$, while Theorem~\ref{thm:root} and $M_n=(1-\delta)n+O(1)$ give
\[
 \eps_{M_n}(\beta;1)
 =\exp\{-2\ell(1-\delta)n+o(n)\}.
\]
Since $\xi<2\ell\delta$, this lower-degree asymptotic is exponentially larger than the upper bound obtained from \eqref{eq:stripped-error}, a contradiction.  Therefore fewer than $\delta n$ poles have height below $s_1$ for all sufficiently large $n$. Since $\delta>0$ was arbitrary, \eqref{eq:boundary-tight} follows.
\end{proof}

\begin{corollary}[Original scaled variables]\label{cor:original-tight}
For every $\epsilon>0$ there exist $0<r<R<\infty$ and an angular separation $\tau>0$, depending only on $\epsilon,\beta$, such that, for all sufficiently large $n$ (also depending on the sequence), at least $(1-\epsilon)n$ of the scaled poles
\[
 \zeta_{j,n}=\frac{hp_{j,n}}n
\]
satisfy
\[
 r\le|\zeta_{j,n}|\le R,
 \qquad
 \dist(\arg\zeta_{j,n},[-\beta/2,\beta/2])\ge\tau
\]
with angles interpreted in the sector complement.  Thus no positive asymptotic pole mass occurs at scales $o(n/h)$ or $\gg n/h$, or on the boundary.
\end{corollary}

\begin{proof}
The inverse of \eqref{eq:Mobius} is a homeomorphism between compact subsets of $\Omega_\beta$ bounded away from $\partial K_\beta$ and compact subsets of the complement of the scaled sector bounded away from its boundary, zero, and infinity.  Apply Theorem~\ref{thm:boundary}.
\end{proof}

\subsection{The Hermite residual identity: turning approximation error into an integral}

After boundary tightness, any subsequential pole limit is an interior probability measure. To identify it, we compare two estimates for a single integral: Green potentials will force a large lower bound if the pole limit is wrong, while the small approximation error will force a smaller upper bound. The following exact Hermite identity connects the two. The zero-table polynomial in the test multiplier is chosen separately from the approximant denominator; the constant term accommodates M\"obius changes of the Cauchy kernel.

\begin{lemma}[Hermite--Walsh residual identity]\label{lem:hermite-residual}
Let $E$ and $F$ be disjoint compact sets and let $\mathcal C$ be a positively oriented Jordan curve with $E$ inside and $F$ outside.  Suppose
\[
 \mathcal F_n(z)=c_n+\frac{1}{2\pi i}\int_F
 \frac{\rho_n(t)\chi(t)}{t-z}\,dt
\]
with $c_n\in\C$ and $\rho_n\chi$ integrable on $F$. Let $R_n=P_n/Q_n$, with both polynomial degrees at most $d$.  Let $\omega_N$ be a polynomial of degree $N$ whose zeros lie on the $E$-side of $\mathcal C$.  If $s$ is a polynomial satisfying
\[
 N\ge d+\deg s+2,
\]
then
\begin{equation}\label{eq:hermite-residual}
 \frac{1}{2\pi i}\int_F Q_n(t)s(t)
 \frac{\rho_n(t)\chi(t)}{\omega_N(t)}\,dt
 =
 \frac{1}{2\pi i}\oint_\mathcal C
 \frac{Q_n(z)(\mathcal F_n(z)-R_n(z))s(z)}{\omega_N(z)}\,dz .
\end{equation}
\end{lemma}

\begin{proof}
The constant term contributes $\oint_\mathcal C c_nQ_ns/\omega_N=0$ by the degree inequality.  Insert the remaining Cauchy representation into the right-hand side and interchange the two integrals.  For $t\in F$ the point $t$ lies outside $\mathcal C$, and the residue theorem on the Riemann sphere gives
\[
 \frac{1}{2\pi i}\oint_\mathcal C
 \frac{Q_n(z)s(z)}{\omega_N(z)(t-z)}\,dz
 =\frac{Q_n(t)s(t)}{\omega_N(t)}.
\]
Indeed, the residue at $z=t$ of the same rational differential is $-Q_n(t)s(t)/\omega_N(t)$ and the residue at infinity vanishes because $N\ge d+\deg s+2$.  Likewise,
\[
 \oint_\mathcal C \frac{P_n(z)s(z)}{\omega_N(z)}\,dz=0,
\]
since $\deg(P_ns)\le N-2$.  Subtraction gives \eqref{eq:hermite-residual}.
\end{proof}

In a canonical multipoint Pad\'e construction, the right-hand side of \eqref{eq:hermite-residual} vanishes, giving the usual varying non-Hermitian orthogonality on the left. Here $R_n$ is arbitrary and no interpolation condition is imposed. The right-hand side need not vanish, but we can bound it using the approximation error and the rational multiplier that we choose.

\subsection{The adaptive Green test}
\label{subsec:adaptive-green}

From now on, $R_n=P_n/Q_n$ denotes the approximant in the current M\"obius chart. The auxiliary $p_n$ and $q_n$ introduced below serve different purposes: $p_n$ supplies mirror zeros, while $q_n$ corrects the local phase (these are not the original numerator and denominator in the sector variable). The table fixes notation; below, local conformal coordinates are $z=\psi(w)$.
\begin{center}
\begin{tabular}{@{}p{.23\linewidth}p{.71\linewidth}@{}}
\toprule
Notation & Meaning\\
\midrule
$E,\ D$ & Compact norming set and its complement $D=\Chat\setminus E$.\\
$F,\ \Gamma$ & Fixed truncated Cauchy contour and its charged arc $\Gamma\subset F$.\\
$\lambda_F,\ \mu$ & Equilibrium probability charge and a candidate pole limit.\\
$\varphi,\ \ell_\beta$ & Transported field $\Rea\Phi_{\rm sec}/2$ and equilibrium level $\ell$.\\
$\sigma,\ \tau$ & Auxiliary zero measure with $|\sigma|<1$, and $\tau=\mu+\sigma$.\\
$Q_n,p_n,q_n$ & Approximant denominator, mirror polynomial, and phase numerator.\\
$\omega_N,\ \chi$ & Boundary-zero interpolation table and fixed Cauchy density.\\
\bottomrule
\end{tabular}
\end{center}
Write $G_D^\tau(z)=\int g_D(z,t)\,d\tau(t)$, $I_D(\tau)=\int G_D^\tau\,d\tau$, and $|\tau|=\tau(D)$. The symbols $\mathsf Z(P)$ and $\capc$ denote the zero counting measure of a polynomial and logarithmic capacity; q.e. means outside a polar set. Standard potential-theoretic facts are used in the conventions of \cite{Ransford1995,Helms2014,SaffTotik1997}.

The Hermite identity is exact, and we are free to choose its zero table. Sweeping the zeros of $Q_np_n$ onto $\partial E$ (balayage) will make the logarithm of the rational multiplier approach $-G_D^\tau\le0$. We must first turn an incorrect pole limit into a strict Green minimum, then realize the resulting potential with a correctly normalized finite table and enough degree reserve, and finally eliminate local complex cancellation without changing that fixed table.

Let $\mu$ be a candidate pole limit. We will construct an auxiliary positive zero measure $\sigma$ with $|\sigma|<1$ and write $\tau=\mu+\sigma$. The desired multiplier has magnitude approximately $\exp(-nG_D^\tau)$, hence at most one on the exponential scale. Proposition~\ref{prop:base-table} supplies the precise upper bound. We do not claim locally uniform convergence of logarithmic potentials from weak zero convergence: for example, $z(z-2)^{n-1}$ has limiting zero measure $\delta_2$ but vanishes at zero for every $n$. Such exceptional zeros only strengthen the one-sided upper estimate used here.

Once the adaptive table has been fixed, we denote the left side of \eqref{eq:hermite-residual} by $\mathcal L_n^{\mathrm{ad}}(s)$. Before constructing that table, we derive the strict exponent gap from $|\sigma|<1$.

\begin{lemma}[Quantitative Green mass gap]\label{lem:green-gap}
Let $\Gamma=\supp\lambda_F\Subset D$ be the regular charged arc, and suppose $G_D^{\lambda_F}+\varphi=\ell_\beta$ on $\Gamma$.  Let $\kappa_\Gamma$ be its unweighted Green equilibrium probability measure, with level $c_\Gamma=I_D(\kappa_\Gamma)>0$.  For every finite positive measure $\tau$ in $D$ with $|\tau|<2$, and every contour $F\supset\Gamma$,
\begin{equation}\label{eq:green-gap}
 \inf_F\{G_D^\tau+2\varphi\}
 \le 2\ell_\beta-(2-|\tau|)c_\Gamma<2\ell_\beta.
\end{equation}
In particular, $\tau$ need not be supported on $F$ or on a compact subset of $D$.
\end{lemma}

\begin{proof}
Write $h_\Gamma=G_D^{\kappa_\Gamma}$.  The unweighted equilibrium principle gives $h_\Gamma=c_\Gamma$ on $\Gamma$ and $0\le h_\Gamma\le c_\Gamma$ in $D$.  By symmetry and Tonelli's theorem,
\[
 \int_\Gamma\varphi\,d\kappa_\Gamma
    =\ell_\beta-\int h_\Gamma\,d\lambda_F=\ell_\beta-c_\Gamma,
 \qquad
 \int_\Gamma G_D^\tau\,d\kappa_\Gamma=\int_D h_\Gamma\,d\tau\le |\tau|c_\Gamma.
\]
Average $G_D^\tau+2\varphi$ over the probability measure $\kappa_\Gamma$ to obtain \eqref{eq:green-gap}.
\end{proof}

Only the charged arc $\Gamma$ is averaged: $\tau$ need not lie on the Cauchy contour, and no tail-potential estimate is used.

\subsection{Finding a strict interior contact point}
\label{subsec:external-selector}

We now show that an incorrect pole limit produces the strict minimum required for the Hermite lower bound. First, the Green potential of $\mu\ne\lambda_F$ falls strictly below that of the equilibrium measure somewhere on the nonseparating charged arc. A small positive perturbation selects one regular point $z_0$ as the unique global minimizer. The following subsection modifies this potential to obtain the local reflection and mirror zeros needed for phase correction.

For this potential-theoretic step only, use a conformal map from $D$ to a bounded domain and push forward the measures; we continue to denote the images by the same symbols. Green functions and total masses are preserved, and the field is transported by composition. No rational function is transformed.  Let $\Gamma$ be the compact analytic simple S-arc.  In particular, $D\setminus\Gamma$ is connected. The equilibrium density already obtained in Section~\ref{sec:Scurve} makes $G_D^{\lambda_F}$ continuous on $\overline D$, with zero boundary values. Let $\kappa_\Gamma$ be the \emph{unweighted} Green equilibrium probability measure of $\Gamma$, and let $c_\Gamma=I_D(\kappa_\Gamma)>0$.  The usual equilibrium and maximum principles give
\begin{equation}\label{eq:unweighted-green-test}
 h_\Gamma:=G_D^{\kappa_\Gamma}=c_\Gamma\quad\hbox{on }\Gamma,
 \qquad 0<h_\Gamma<c_\Gamma\quad\hbox{in }D\setminus\Gamma.
\end{equation}
Regularity removes the exceptional boundary set, and $\kappa_\Gamma$ has full support on $\Gamma$; see also \cite[Sec.~2]{Pouliasis2011}. In particular, the strict inequality away from the nonseparating arc lets us detect an incorrect pole measure even when it has mass off the arc.

\begin{lemma}[Negative contact on a nonseparating support]
\label{lem:negative-contact}
Let $\lambda_F$ be any probability measure on $\Gamma$ with continuous Green potential, and let $\mu$ be any probability measure in $D$.  No compact-support or finite-energy assumption on $\mu$ is required.  If $\mu\ne\lambda_F$, then
\begin{equation}\label{eq:negative-contact}
 \capc\{z\in\Gamma:G_D^\mu(z)-G_D^{\lambda_F}(z)<0\}>0.
\end{equation}
Moreover, for $W_0=G_D^\mu-G_D^{\lambda_F}$,
\begin{equation}\label{eq:negative-global-min}
 m_0:=\min_\Gamma W_0<0,\qquad W_0(z)>m_0\quad(z\in D\setminus\Gamma).
\end{equation}
The infimum over the regular interior $\Gamma^\circ$ is also $m_0$.
\end{lemma}

\begin{proof}
Tonelli's theorem and \eqref{eq:unweighted-green-test} give the finite identity
\begin{equation}\label{eq:negative-average}
 \int_\Gamma W_0\,d\kappa_\Gamma
   =\int_D(h_\Gamma-c_\Gamma)\,d\mu\le0.
\end{equation}
If $\mu$ puts positive mass outside $\Gamma$, the right side is strictly negative. Since $\kappa_\Gamma$ has finite energy and does not charge polar sets, the negative set on $\Gamma$ must have positive capacity, proving \eqref{eq:negative-contact} in this case.

The remaining case is $\supp\mu\subset\Gamma$, for which the average in \eqref{eq:negative-average} vanishes. Suppose, toward a contradiction, that $W_0\ge0$ quasi-everywhere on $\Gamma$.  The integral in \eqref{eq:negative-average} is zero, so $W_0=0$ almost everywhere for $\kappa_\Gamma$.  Lower semicontinuity of $W_0$ and full support of $\kappa_\Gamma$ exclude $W_0>0$ at any point of $\Gamma$.  Therefore $G_D^\mu\le G_D^{\lambda_F}$ everywhere on $\Gamma$, and $\mu$ has finite energy. Both measures now ignore polar sets, and
\[
 I_D(\mu-\lambda_F)=\int_\Gamma W_0\,d(\mu-\lambda_F)=0.
\]
Strict positive definiteness of the Green energy gives $\mu=\lambda_F$, a contradiction.  This proves \eqref{eq:negative-contact} also in the remaining case.

Having established negative contact, we locate its global minimum. The function $W_0$ is lower semicontinuous and bounded below by the negative of the bounded continuous potential $G_D^{\lambda_F}$; hence it attains a negative minimum on the compact arc $\Gamma$.  On $D\setminus\Gamma$ it is superharmonic, its boundary lower limit at $\Gamma$ is at least $m_0$, and at $\partial D$ its lower limit is nonnegative, since $G_D^\mu\ge0$ and $G_D^{\lambda_F}\to0$.  The minimum principle, followed by the strong minimum principle in the connected slit domain, gives \eqref{eq:negative-global-min}.

We also need regular interior points arbitrarily close to this minimum; a minimum confined to an endpoint would not suffice for local reflection.  Suppose to the contrary that $m_0<\inf_{\Gamma^\circ}W_0$, and choose $m_0<b<\min\{0,\inf_{\Gamma^\circ}W_0\}$. The extended minimum principle, ignoring the two polar endpoints, gives $W_0\ge b$ in $D\setminus\Gamma$.  At an endpoint attaining $m_0$, the supermean inequality for $G_D^\mu$ on shrinking circles, together with continuity of $G_D^{\lambda_F}$, would give $W_0\ge b$.  This is a contradiction.
\end{proof}

\begin{remark}[Why nonseparation is essential]
To see why the geometry matters, take $D=\{|z|<1\}$, let $\lambda_F$ be uniform arclength on $|z|=r<1$, and let $\mu=\delta_0$.  Direct angular averaging gives
\[
 G_D^{\lambda_F}(z)=-\log\max\{r,|z|\},\qquad
 G_D^\mu(z)=-\log|z|.
\]
The potentials agree outside the circle and on it, and $\mu$ sweeps to $\lambda_F$. The negative-contact conclusion fails because a closed circle separates $D$. The S-support used here is a compact simple arc, not such a circle.
\end{remark}

\begin{lemma}[Regular strict-minimum Green selector]
\label{lem:green-selector}
Let $F\Subset D$ be a compact simple contour containing $\Gamma$, and let $\varphi$ be harmonic on $D$.  Suppose $\lambda_F$ has the continuity properties above.  For every probability measure $\mu\ne\lambda_F$ in $D$, there exist $0<t<1/2$, a compactly supported probability measure $\sigma_0$, and $z_0\in\Gamma^\circ$ such that
\begin{equation}\label{eq:green-selector-W}
 W=G_D^{\mu+t\sigma_0}+(1+2t)\varphi
       -(1-2t)(G_D^{\lambda_F}+\varphi)
   =W_0+t(2G_D^{\lambda_F}+4\varphi+G_D^{\sigma_0})
\end{equation}
has a unique finite global minimum on $F$ at $z_0$.  It also has a unique minimum there on a closed Jordan neighborhood of $\Gamma$.  The point $z_0$ is not an atom of $\mu$, and $\sigma_0$ has no support in a neighborhood of $z_0$.
\end{lemma}

\begin{proof}
Choose a smooth Jordan domain $K$ with $\Gamma\Subset K\Subset D$. By Lemma~\ref{lem:negative-contact} and compactness, the minimum $m_0$ on $\Gamma$ lies strictly below the minimum of $W_0$ on $A=\partial K\cup(F\setminus K)$. We will preserve this gap while perturbing $W_0$. Let the endpoints of $\Gamma$ be $a_\pm$, and put $\sigma_{\mathrm{end}}=(\delta_{a_+}+\delta_{a_-})/2$.  Choose an open disk $B$ with $\overline B\Subset D\setminus(\overline K\cup F)$, a point $b_0\in B$, and a dense sequence $b_j\in B$.  For independent parameters $0\le\alpha_j\le2^{-j-2}$ define
\[
 \sigma_{\mathrm{ext}}=\left(1-\sum_{j\ge1}\alpha_j\right)\delta_{b_0}
             +\sum_{j\ge1}\alpha_j\delta_{b_j},
 \qquad \sigma_0=\tfrac12\sigma_{\mathrm{end}}+\tfrac12\sigma_{\mathrm{ext}}.
\]
These are positive probability measures, and $G_D^{\sigma_{\mathrm{ext}}}$ is uniformly bounded on $\overline K$, independently of all parameters. The endpoint atoms will exclude endpoint minimizers; the exterior atoms provide small adjustable harmonic perturbations along $\Gamma$.

Choose $z_1\in\Gamma^\circ$ with finite $W_0(z_1)$ sufficiently close to $m_0$; this is possible by the last assertion of Lemma~\ref{lem:negative-contact}. The endpoint potential is finite at $z_1$.  Since $G_D^{\lambda_F},G_D^{\sigma_0}\ge0$ and $\varphi$ is bounded on the relevant compact sets, for all sufficiently small fixed $t>0$, uniformly in the parameters,
\[
 W(z_1)<\inf_A W.
\]
Since $W$ is superharmonic on $K\setminus\Gamma$ and lower semicontinuous on $\overline K$, the minimum principle places its minimum on $\Gamma$. The gap against $A$ makes this also the global minimum on $F$. The endpoint atoms give $W(a_\pm)=+\infty$, so every minimizer lies in the regular interior $\Gamma^\circ$.

It remains to distinguish a single minimizer, without making any nondegeneracy assumption on $W_0$. We use the exterior point charges for this purpose.  On $\Gamma$ the parameter dependence is
\[
 G_D^{\sigma_{\mathrm{ext}}}(x)=g_D(x,b_0)+\sum_{j\ge1}\alpha_j
 \{g_D(x,b_j)-g_D(x,b_0)\}.
\]
Because $\overline B$ and $\Gamma$ are disjoint compact sets, the functions
\[
 f_j(x)=g_D(x,b_j)-g_D(x,b_0),\qquad x\in\Gamma,
\]
are uniformly bounded.  They also separate points of $\Gamma$.  Indeed, equality $f_j(x)=f_j(y)$ for every $j$, followed by density of $b_j$ in $B$, would make $g_D(x,b)-g_D(y,b)$ constant on $B$, hence, by harmonic unique continuation, on $D\setminus\{x,y\}$.  Its logarithmic singularities make this impossible unless $x=y$.

The terms of $W|_\Gamma$ independent of the $\alpha_j$ form a lower semicontinuous base function with a finite minimum; the endpoint atoms assign $+\infty$ at $a_\pm$.  Since $\sum_j2^{-j-2}<\infty$, Lemma~\ref{lem:generic-unique}, applied to the perturbations $(t/2)f_j$, gives one parameter vector for which $W$ has exactly one minimizer $z_0$ on $\Gamma$.  The strong minimum principle excludes any other minimizing point in $K\setminus\Gamma$.

Finiteness at $z_0$ excludes an atom of $\mu$ there.  Finally, $\supp\sigma_0\subset\{a_-,a_+\}\cup\overline B$, so a neighborhood of $z_0$ is disjoint from $\supp\sigma_0$.
\end{proof}

\subsection{Green gluing and local reflection}
\label{subsec:external-gluing}

The selector gives a unique interior minimum, but not yet the mirror symmetry required for the phase argument. Following Gonchar--Rakhmanov \cite{GoncharRakhmanov1989}, we modify the potential only near that minimum. We state the Green version explicitly because the candidate measure may lack compact support and finite energy. Two preliminary lemmas explain what survives this modification: total Riesz mass and local reflection of the measure.

\begin{lemma}[Mass conservation under compact Green gluing]
\label{lem:green-compact-mass}
Let $\alpha$ be a finite positive measure with compact support in the bounded regular domain $D$.  If $u$ is superharmonic in $D$ and $u-G_D^\alpha$ has compact support in $D$, then
\[
 u=G_D^\xi,\qquad \xi\ge0,\qquad |\xi|=|\alpha|.
\]
\end{lemma}

\begin{proof}
Put $v=u-G_D^\alpha$ and $\xi=-(2\pi)^{-1}\Delta u$. Since $u$ is superharmonic and agrees with $G_D^\alpha$ outside a compact subset of $D$, its Riesz measure $\xi$ is positive and compactly supported.  Distributionally, $\xi-\alpha=-(2\pi)^{-1}\Delta v$.  Testing against a smooth cutoff equal to one near $\supp v\cup\supp\alpha$ gives $|\xi|-|\alpha|=0$. The minimum principle gives $u\ge0$, since it equals the nonnegative Green potential outside a compact set.  Its Riesz decomposition is $G_D^\xi+H$ with $H$ nonnegative harmonic.  Near the boundary $H\le u=G_D^\alpha\to0$ uniformly; hence $H=0$.
\end{proof}

\begin{lemma}[Local S-reflection and reflected Riesz measure]
\label{lem:local-reflection}
Let $z_0$ be a regular point of a real-analytic S-arc $\Gamma\Subset D$.  Let $\mathcal T$ be continuous across $\Gamma$, harmonic on the two sides, constant on $\Gamma$, and suppose
\[
 \partial_{n_+}\mathcal T=\partial_{n_-}\mathcal T\quad\hbox{on }\Gamma,
\]
where $n_+$ and $n_-$ point from the arc into its two opposite sides, rather than denoting the same oriented normal on both sides.  Then, after shrinking a neighborhood of $z_0$, there is an anticonformal involution $z\mapsto z^*$ fixing $\Gamma$ pointwise such that
\begin{equation}\label{eq:T-reflection}
 \mathcal T(z^*)=\mathcal T(z).
\end{equation}
Moreover, if $\alpha$ is a finite positive measure and $V$ is a sufficiently small reflection-invariant neighborhood of $z_0$, then
\begin{equation}\label{eq:green-reflection-harmonic}
 z\longmapsto G_D^{\alpha|_V}(z^*)-G_D^{(\alpha|_V)^*}(z)
\end{equation}
admits a harmonic extension to a smaller neighborhood of $z_0$.  Here the difference is initially understood where both Green potentials are finite; at their common infinite points it means the extension obtained by cancelling the singular kernels before integration, not a formal $+\infty-(+\infty)$.  Consequently the local Riesz measure of $G_D^\alpha(z^*)$ is the reflected measure $\alpha^*$.
\end{lemma}

\begin{proof}
Choose a conformal coordinate $z=\psi(w)$ that maps a real interval to $\Gamma$, with $\psi(0)=z_0$, and define $z^*=\psi(\overline w)$.  Write $U_\pm=\mathcal T\circ\psi$ on the upper and lower half-disks.  Both have the same constant boundary value on the real interval.  Equality of these side-normal derivatives becomes
\[
 \partial_yU_+(x,0)=-\partial_yU_-(x,0).
\]
After subtracting the constant boundary value, odd harmonic reflection extends each branch across the real interval.  The analytic harmonic functions $U_+(x,y)$ and $U_-(x,-y)$ have the same Cauchy data there, so local uniqueness gives $U_+(x,y)=U_-(x,-y)$, which is \eqref{eq:T-reflection}.

To reflect an arbitrary measure, we compare the Green kernels rather than subtracting two potentially infinite values. In the same coordinate, for $w,\xi$ in a smaller disk,
\[
 g_D(\psi(w),\psi(\xi))=-\log|w-\xi|+h(w,\xi),
\]
where $h$ is jointly continuous and bounded on smaller closed product disks, and harmonic in its first variable. Therefore
\[
 g_D(\psi(\overline w),\psi(\xi))
 -g_D(\psi(w),\psi(\overline\xi))
\]
is harmonic in $w$: the singular terms cancel because $|\overline w-\xi|=|w-\overline\xi|$, and the two regular parts are harmonic (the first after anticonformal composition).  The logarithmic terms cancel \emph{before} integration against the pull-back of $\alpha|_V$: the result is the integral of the two bounded regular kernels, hence a finite harmonic function on the smaller disk.  Where the individual potentials are finite it equals the displayed difference.  This defines its harmonic extension also where both potentials are $+\infty$, without assigning a value to their formal difference.  Taking local Riesz measures then gives the last assertion.
\end{proof}

We now combine the unique minimizer with local reflection. As in \cite[\S3, Lemma~9]{GoncharRakhmanov1989}, taking the minimum of two superharmonic functions lets us replace a potential near $z_0$ without changing it far away. The Green formulation also tracks positivity, total mass, and the strict reserve $|\sigma|<1$: this reserve will later control the degree of the mirror polynomial and yield the strict exponent gap.

\begin{theorem}[Green adaptation of the GR reflection construction]
\label{thm:green-GL9}
Let $D,F,\Gamma,\lambda_F,\varphi$ be as in Lemma~\ref{lem:green-selector}, and suppose
\[
 \mathcal T=G_D^{\lambda_F}+\varphi=\ell_\beta\quad\hbox{on }\Gamma,
 \qquad \mathcal T\ge\ell_\beta\quad\hbox{on }F,
\]
with the S-property on the analytic arc $\Gamma$.  Then the local reflection $z\mapsto z^*$ is supplied by Lemma~\ref{lem:local-reflection}.  For every probability measure $\mu\ne\lambda_F$ in $D$, there exist a compactly supported positive measure $\sigma$, $|\sigma|<1$, and a regular point $z_0\in\Gamma^\circ$ such that
\[
 H=G_D^{\mu+\sigma}+2\varphi
\]
has a unique finite global minimum on $F$ at $z_0$.  In a reflection-invariant neighborhood $U$ of $z_0$,
\begin{equation}\label{eq:GL9-local-complete}
 H(z^*)=H(z),\qquad \sigma|_U=(\mu|_U)^*.
\end{equation}
Theorem~\ref{thm:adaptive-sync} realizes this measure by mirror polynomials. The strict mass bound will give the exponent gap in Lemma~\ref{lem:green-gap}.
\end{theorem}

\begin{proof}
We build the reflected auxiliary measure in three steps. The function $W_1$ contains the selected minimum. A second function $W_2$ contains the locally reflected part of $\mu$ and a slightly reduced equilibrium background; its missing mass provides the degree reserve. We then glue the reflected $W_1$ to $W_2$ by taking a minimum, and recover a positive measure of controlled total mass from Lemma~\ref{lem:green-compact-mass}.

Take $t,\sigma_0,z_0,W$ from Lemma~\ref{lem:green-selector}, and set
\[
 W_1=G_D^{\mu+t\sigma_0}+(1+2t)\varphi,
 \qquad W=W_1-(1-2t)\mathcal T.
\]
Choose a small reflection-invariant neighborhood $U_1$ of $z_0$, compactly contained in the selector neighborhood, with $F\cap U_1\subset\Gamma$ and $\supp\sigma_0\cap\overline U_1=\varnothing$.  By strict minimality and lower semicontinuity,
\[
 d:=\inf_{\partial U_1}(W-W(z_0))>0.
\]
For a much smaller reflection-invariant neighborhood $V\Subset U_1$, put $\nu=(\mu|_V)^*$.  Since $z_0$ is not an atom, $V$ may be chosen so that
\begin{equation}\label{eq:GL9-nu-control}
 |\nu|<t/2,\qquad \sup_{\partial U_1}G_D^\nu<d/4.
\end{equation}
The second inequality holds because the Green kernel is bounded when the source approaches $z_0$ but the evaluation point stays on the fixed outer boundary. Define
\[
 C=W(z_0)+d/2,\qquad W_2=(1-2t)\mathcal T+G_D^\nu+C.
\]
We first compare the two proposed branches near the outer boundary of $U_1$. Where their values are finite, reflection symmetry of $\mathcal T$ gives
\[
 W_1(z^*)-W_2(z)=W(z^*)-W(z_0)-G_D^\nu(z)-d/2.
\]
At common infinite points this is understood through the regular-part comparison below, not as a formal difference. It is positive near $\partial U_1$ by \eqref{eq:GL9-nu-control}.  Near $z_0$ the comparison must be made \emph{after} cancelling the common singularity.  On a smaller invariant disk $V_0\Subset V$, write
\[
 P(z)=G_D^\nu(z),\qquad
 W_1(z^*)=P(z)+h_1(z),\qquad W_2(z)=P(z)+h_2(z).
\]
By Lemma~\ref{lem:local-reflection}, $h_1$ is continuous: the local reflected $\mu$-potential minus $P$ has a harmonic extension, and all other source terms are supported away from $V_0$.  Also $h_2=(1-2t)\mathcal T+C$ is continuous.  Finiteness of $G_D^\mu(z_0)$ and local bi-Lipschitz bounds for reflection give $P(z_0)<\infty$ and
\[
 h_1(z_0)-h_2(z_0)=-P(z_0)-d/2<0.
\]
After shrinking $V_0$, we have $h_1<h_2$ throughout $V_0$. This determines the $W_1(z^*)$ branch by comparison of the finite regular parts, even where $P=+\infty$; no strict inequality between two infinite potentials is asserted. The branch choice is therefore meaningful without continuity or finite energy of $G_D^\mu$.

Choose a still smaller invariant $U_0\Subset V$ and glue
\[
 W_3(z)=
 \begin{cases}
  W_1(z^*),&z\in U_0,\\
  \min\{W_1(z^*),W_2(z)\},&z\in U_1\setminus\overline U_0,\\
  W_2(z),&z\in D\setminus U_1.
 \end{cases}
\]
The inner and outer comparisons ensure that these pieces agree in the required order at the gluing interfaces. The minimum of superharmonic functions is superharmonic, so $W_3$ is superharmonic. Since $\varphi$ is harmonic, $u=W_3-(1-2t)\varphi-C$ is superharmonic and agrees outside $U_1$ with $(1-2t)G_D^{\lambda_F}+G_D^\nu$.  Lemma~\ref{lem:green-compact-mass} gives
\[
 u=G_D^{\sigma_1},\qquad \sigma_1\ge0,
 \qquad |\sigma_1|=1-2t+|\nu|.
\]
Set $\sigma=\sigma_1+t\sigma_0$.  Then
\begin{equation}\label{eq:GL9-mass-complete}
 |\sigma|=1-t+|\nu|<1-t/2<1,
 \qquad H=W_1+W_3-C.
\end{equation}
On $U_0$ we have $W_3(z)=W_1(z^*)$, so $H(z)=W_1(z)+W_1(z^*)-C$ is reflection-symmetric. Taking local Riesz measures and using the absence of $\sigma_0$-mass near $z_0$ gives $\sigma|_U=(\mu|_U)^*$ after shrinking $U$.

At $z_0$, $H(z_0)=2W_1(z_0)-C$.  On the part of $F$ where $W_3=W_2$,
\[
 H(z)-H(z_0)
 =[W(z)-W(z_0)]+2(1-2t)[\mathcal T(z)-\ell_\beta]
       +G_D^\nu(z)+d/2>0.
\]
On the other branch, $z\in\Gamma$ and $z^*=z$, so the difference is $2[W(z)-W(z_0)]>0$ for $z\ne z_0$.  Thus the gluing leaves $z_0$ as the unique global minimizer of $H$ on $F$, while enforcing the desired local reflection.
\end{proof}

\begin{corollary}[Mirror polynomials]
\label{cor:mirror-degree}
Suppose $n^{-1}\mathsf Z(Q_n)\weakto\mu$, and let $\sigma,U$ be supplied by Theorem~\ref{thm:green-GL9}.  The neighborhood can be chosen with $\mu(\partial U)=\sigma(\partial U)=0$.  There exist polynomials $p_n$ satisfying
\[
 n^{-1}\mathsf Z(p_n)\weakto\sigma,\qquad
 \deg p_n/n\longrightarrow|\sigma|<1,
\]
whose zeros in $U$ are exactly the reflected zeros of $Q_n$ there, with multiplicity.
\end{corollary}

\begin{proof}
First include, with multiplicity, the reflected zeros of $Q_n$ lying in $U$. Their normalized counting measures converge to $(\mu|_U)^*=\sigma|_U$ because $U$ is a continuity set for both measures. Then discretize the remaining positive measure $\sigma|_{D\setminus U}$ using zeros outside $U$.  Rounding changes the degree by $o(n)$, and the total normalized degree tends to $|\sigma|$.  Thus for some fixed $\delta>0$, $\deg p_n\le(1-\delta)n$ for all sufficiently large $n$.
\end{proof}

\subsection{Building the adaptive table with a controlled normalization}
\label{subsec:adaptive-sync}

The Green construction gives measures, whereas the Hermite identity requires actual polynomials. We therefore discretize the auxiliary and swept measures and fix the scalar normalization of the resulting quotient. Zero locations alone determine a logarithmic potential only up to an additive constant; that constant matters at exponential accuracy. Both the base table and its normalization are fixed before the phase correction; that correction only appends its own prescribed factor to the table.

From this point on, we work in a finite affine coordinate selected by a fixed M\"obius transformation; no general conformal map is applied to the rational functions.  For a given subsequence and its limiting measure, its point at infinity can be chosen in $D\setminus(F\cup\supp\sigma)$, outside the countable set of poles of the approximants and the atoms of the limiting measure.  Then $E,F$ are compact in this coordinate, $\deg Q_n=d_n=\operatorname{ord}(R_n)$, and $\deg P_n\le d_n$.  The bounded conformal coordinate in the proof of Theorem~\ref{thm:green-GL9} was used only for potentials and did not change the rational approximation class.  A M\"obius-transformed finite Cauchy transform has the same form with an additional constant, already allowed in Lemma~\ref{lem:hermite-residual}.

We begin by fixing the denominator normalization.  If the finite zeros of $Q_n$ are $a_{j,n}$, rescale $P_n,Q_n$ together so that
\begin{equation}\label{eq:Q-chordal-normalization}
 Q_n(z)=\prod_{j=1}^{d_n}\frac{z-a_{j,n}}{(1+|a_{j,n}|^2)^{1/2}}.
\end{equation}
This joint rescaling leaves $R_n$ unchanged but fixes the scalar of $Q_n$. We use the same initial normalization for the auxiliary zero products. Proposition~\ref{prop:base-table} will then apply one explicit exponential scalar to $p_n$ to correct the balayage constant. The later phase correction does not change this choice.

Using chordal distance $[z,a]$, set $V_*^\tau(z)=\int\log(1/[z,a])\,d\tau(a)$, and put $L(a)=\tfrac12\log(1+|a|^2)$ for finite $a$.  For $a\in D$, let $\omega_a$ be harmonic measure on $\partial E$ viewed from $a$, with the point at infinity included as an ordinary point of the spherical domain $D$. The logarithmic sweep used below is mass preserving: $\tau^E=\int\omega_a\,d\tau(a)$.  It is different from sweeping onto an interior contour, where mass may be lost at the fixed boundary.

\begin{lemma}[Bounded spherical balayage constant]
\label{lem:spherical-balayage}
Suppose $E$ is a regular compact continuum with connected complement $D$. For finite $a\in D$, define
\[
 c_E(a)=L(a)-g_D(a,\infty)-\int_{\partial E}L(t)\,d\omega_a(t).
\]
This function extends continuously and boundedly to $\overline D$ on the sphere, with $c_E=0$ on $\partial E$ and
\[
 c_E(\infty)=\log\capc(E)-\int L(t)\,d\omega_\infty(t).
\]
For every finite positive measure $\tau$ in $D$, with no moment or energy assumption,
\begin{equation}\label{eq:spherical-balayage}
 V_*^\tau-V_*^{\tau^E}=G_D^\tau+C_E(\tau),
 \qquad C_E(\tau)=\int c_E(a)\,d\tau(a)\in\R.
\end{equation}
\end{lemma}

For a point charge this is the usual Green representation with the chordal normalizing terms included. Its important feature is the boundedness of $c_E$ even at infinity: we can integrate the constant against measures with no logarithmic-moment assumption. The proof is in \ref{app:balayage-constant}.

\begin{proposition}[A base table fixed before phase correction]
\label{prop:base-table}
Let $n^{-1}\mathsf Z(Q_n)\weakto\mu$ with $d_n/n\to1$, and choose $p_n$ as in Corollary~\ref{cor:mirror-degree}.  Write $\tau=\mu+\sigma$ and $d_n^{\mathrm{tot}}=d_n+\deg p_n$.  There is a polynomial $\omega_n^0$ with all zeros on $\partial E$ and degree $N_n^0=d_n^{\mathrm{tot}}+2$, and a scalar normalization of $p_n$, such that, for $A_n=Q_np_n/\omega_n^0$,
\begin{equation}\label{eq:base-table-limits}
 \begin{split}
 n^{-1}\mathsf Z(\omega_n^0)&\weakto\tau^E,\\
 n^{-1}\log|A_n|&\longrightarrow-G_D^\tau
                         \quad\hbox{in }L^1_{\mathrm{loc}}(D),\\
 \limsup_{n\to\infty}n^{-1}\log\|A_n\|_K
       &\le-\inf_KG_D^\tau\le0\qquad(K\Subset D).
 \end{split}
\end{equation}
The last two statements are understood in finite coordinate neighborhoods. In particular, for $z_n\to z\in D$ in such a neighborhood,
\begin{equation}\label{eq:base-moving-envelope}
 \limsup_n n^{-1}\log|A_n(z_n)|\le-G_D^\tau(z).
\end{equation}
The table depends on $Q_n,p_n$, but not on a later phase corrector.
\end{proposition}

\begin{proof}
The construction has two independent parts. We discretize the swept zero measure, adding exactly two boundary zeros to meet the Hermite degree condition. Then we fix the exponential scalar of $p_n$ so that the resulting logarithmic exponent is $-G_D^\tau$, not $-G_D^\tau$ plus an unknown constant.

Put $\tau_n=n^{-1}\mathsf Z(Q_np_n)$, so $|\tau_n|=d_n^{\mathrm{tot}}/n$ and $\tau_n\weakto\tau$.  Quantize the probability measure $\tau_n^E/|\tau_n|$ by $d_n^{\mathrm{tot}}$ equal-weight points on the compact set $\partial E$, with weak discrepancy tending to zero, and append two copies of a fixed $a\in\partial E$.  Such quantization can be obtained by partitioning $\partial E$ into sets of diameter tending to zero, using $o(d_n^{\mathrm{tot}})$ cells, and rounding the cell counts.  Weak continuity of harmonic sweeping then gives the first limit in \eqref{eq:base-table-limits}.

Use the chordal zero-product normalization \eqref{eq:Q-chordal-normalization} for $p_n$ before its final scalar and for $\omega_n^0$.  Multiply the normalized auxiliary numerator product by $\exp\{nC_E(\tau_n)\}$, which changes only the scalar of $p_n$ and is the single exponential renormalization allowed by the convention above.  With $\nu_n=n^{-1}\mathsf Z(\omega_n^0)$ one has the exact formula
\begin{equation}\label{eq:base-normalization-exact}
 \frac1n\log|A_n(z)|
 =-V_*^{\tau_n}(z)+V_*^{\nu_n}(z)+C_E(\tau_n)-\frac{2}{n}L(z).
\end{equation}
The constants converge by Lemma~\ref{lem:spherical-balayage}. Since the table zeros lie on $\partial E$, its logarithmic potentials converge locally uniformly in $D$; the numerator potentials converge in local planar $L^1$.  For the latter statement, truncate the logarithmic singularity and use that its local area integral is uniformly small with respect to the source point.  Equation~\eqref{eq:spherical-balayage} gives the second limit. Portmanteau's upper-envelope inequality for the bounded-above kernel $\log[z,a]$ gives the compact-norm bound. Applying it on $\overline{B(z,r)}\Subset D$ and then letting $r\downarrow0$ gives \eqref{eq:base-moving-envelope}, since $G_D^\tau$ is lower semicontinuous. Individual exceptional zeros are allowed; local uniform convergence of their logarithmic potentials is not needed.
\end{proof}

\subsection{Removing complex cancellation near the strict minimum}
\label{subsec:rational-phase}

The Green minimum controls the exponential size of the integrand, but the complex integral could still be much smaller through cancellation. The reflection argument of Gonchar--Rakhmanov \cite{GoncharRakhmanov1989} removes this obstacle locally. We obtain a rational phase corrector of degree $o(n)$ with one prescribed pole; inserting its numerator into the test and its denominator into the table preserves the exact Hermite identity.

\begin{lemma}[Logarithmic means along a regular arc]
\label{lem:arc-log-means}
Let $z=\psi(x)$, $|x|\le r_0$, parametrize a regular analytic arc, with $\psi(0)=z_0$.  Chordal logarithmic potentials of weakly convergent finite positive measures on the sphere converge in $L^1([-r_0,r_0],dx)$ after composition with $\psi$.  If $G_D^\tau(z_0)<\infty$, then
\begin{equation}\label{eq:arc-mean-at-point}
 \lim_{r\downarrow0}\frac1{2r}\int_{-r}^r
 G_D^\tau(\psi(x))\,dx=G_D^\tau(z_0).
\end{equation}
Consequently, suppose $H=G_D^\tau+2\varphi$ has a finite minimum $m_*$ on a compact contour $F$ at $z_0$.  Let $I$ be a sufficiently small fixed real interval about zero with $\psi(I)\subset F$.  If nonnegative $f_n$ satisfy
\[
 n^{-1}\log f_n\longrightarrow-H\circ\psi\quad\hbox{in }L^1(I),
 \qquad \limsup n^{-1}\log\|f_n\|_I\le-m_*,
\]
then
\begin{equation}\label{eq:local-laplace}
 \left(\int_I f_n(x)\,dx\right)^{1/n}\longrightarrow e^{-m_*}.
\end{equation}
\end{lemma}

The logarithmic singularity is uniformly integrable along a short regular analytic arc, which gives the $L^1$ convergence. To obtain the lower bound for the integral, average on intervals shrinking to $z_0$ and apply Jensen's inequality. The details are in \ref{app:arc-means}.

The phase-correction idea comes from \cite[\S3, Lemma~8]{GoncharRakhmanov1989}, where a polynomial of degree $o(n)$ removes local oscillation. We use a rational version whose only pole is prescribed on $E$. This makes it possible to add the denominator factor to the Hermite table at exactly the same time as the numerator enters the test polynomial.

\begin{lemma}[Rational local phase correction with a prescribed pole]
\label{lem:rational-phase}
Let a conformal coordinate $z=\psi(w)$ straighten the local reflection to $w\mapsto\overline w$ on a fixed symmetric disk, with $\psi(0)=z_0$. Suppose analytic differentials $J_n(z)\,dz$ have, there, exactly paired zeros:
\begin{equation}\label{eq:paired-factorization}
 J_n(\psi(w))\psi'(w)=D_n(w)B_n(w),\qquad
 D_n(w)=\prod_j(w-w_{j,n})(w-\overline{w_{j,n}}),
\end{equation}
where $B_n$ is analytic and nonvanishing.  Assume $n^{-1}\log|J_n|\to-H$ in local planar $L^1$, where $H(\psi(\overline w))=H(\psi(w))$.  Fix a finite point $a$ outside the local disk. For a sufficiently small fixed segment $J=\psi(I)\subset F$ about $z_0$, where $I$ is a real interval, there exist
\begin{equation}\label{eq:rational-corrector}
 \vartheta_n(z)=\frac{q_n(z)}{(z-a)^{k_n}},\qquad
 \deg q_n\le k_n=o(n),
\end{equation}
with $\|\vartheta_n\|_K\le e^{o(n)}$ for every fixed compact $K\subset\Chat\setminus\{a\}$, such that, in the local real parameter,
\begin{equation}\label{eq:positive-local-integral}
 \Rea\{\vartheta_n(\psi(x))J_n(\psi(x))\psi'(x)\}
 \ge(1-o(1))|J_n(\psi(x))\psi'(x)|\qquad(x\in I).
\end{equation}
The orientation is chosen so that the real parameter increases along $F$.
\end{lemma}

\begin{proof}
First, local reflection gives an analytic corrector $\beta_n$ of modulus one on the straightened real segment; its product with the integrand is nonnegative there. Second, a Taylor polynomial in $v=(z-z_0)/(z-a)$ approximates this corrector. As a function of $z$, that polynomial is rational with the prescribed pole $a$, and its degree is $o(n)$.

Choose an analytic logarithm $b_n=\log B_n$ on the symmetric disk and set $b_n^*(w)=\overline{b_n(\overline w)}$.  Since $\log|D_n(w)|=\log|D_n(\overline w)|$, the harmonic functions
\[
 \frac1n\{\log|B_n(w)|-\log|B_n(\overline w)|\}
\]
tend to zero in local $L^1$, hence locally uniformly.  Define
\[
 \beta_n(w)=\exp\{(b_n^*(w)-b_n(w))/2\}.
\]
Then $\|\beta_n\|\le e^{o(n)}$ on every smaller disk, while for real $x$
\[
 |\beta_n(x)|=1,\qquad
 \beta_n(x)D_n(x)B_n(x)=D_n(x)|B_n(x)|\ge0.
\]
This remains valid at real zeros: they occur with even multiplicity in $D_n$, so we never divide by a factor vanishing on the arc.

Use the M\"obius coordinate $v=(z-z_0)/(z-a)$ and let $f_n(v)=\beta_n(\psi^{-1}(z(v)))$ on a fixed disk $|v|<R_1$. Choose fixed $0<r<R<R_1$ so that the segment $J$ lies in $|v|\le r$. On $|v|=R$, let $M_n=\max|f_n|=e^{o(n)}$; also $M_n\ge1$ since $|f_n(0)|=1$.  If $f_n(v)=\sum_{j\ge0}c_{j,n}v^j$, Cauchy's inequality gives
\[
 |c_{j,n}|\le M_nR^{-j},\qquad
 \left|f_n(v)-\sum_{j=0}^{k}c_{j,n}v^j\right|
 \le\frac{M_n(r/R)^{k+1}}{1-r/R}\quad(|v|\le r).
\]
Choose $k_n$ large enough to make this error at most $1/n$.  Since $\log M_n=o(n)$, one can take $k_n=O(\log M_n+\log(n+1))=o(n)$. Now set
\begin{equation}\label{eq:corrector-polynomial-explicit}
 q_n(z)=\sum_{j=0}^{k_n}c_{j,n}(z-z_0)^j(z-a)^{k_n-j}.
\end{equation}
This gives \eqref{eq:rational-corrector} and uniform approximation to the exact unit-modulus phase on $J$, proving \eqref{eq:positive-local-integral}. For any fixed $K$ avoiding $a$, $|v|$ is bounded, including at infinity, and
\[
 \|\vartheta_n\|_K\le M_n(k_n+1)
       \max\{1,\|v\|_K/R\}^{k_n}=e^{o(n)}.
\]
The same explicitly constructed rational function therefore gives both the local sign correction and the global subexponential bound.
\end{proof}

\subsection{The final synchronized contradiction}
\label{subsec:sync-rigidity}

\noindent\emph{Order of choices.} First fix the truncated Cauchy contour and a subsequence with pole limit $\mu$. Choose $\sigma,z_0$ and the mirror polynomials $p_n$, then the base table $\omega_n^0$ and its scalar normalization. Only afterward choose the phase corrector $q_n/(z-a)^{k_n}$, inserting $(z-a)^{k_n}$ into the final table simultaneously. This fixes one test polynomial, one table, and one Hermite functional. In Proposition~\ref{prop:sync-residual-upper} we may vary the separating contour with a tolerance $\epsilon$, but not the functional.

\begin{theorem}[Synchronized adaptive Hermite test]
\label{thm:adaptive-sync}
Let $\mu\ne\lambda_F$ be an interior subsequential pole limit, and let $\sigma,z_0$ be furnished by Theorem~\ref{thm:green-GL9}.  Set
\[
 \tau=\mu+\sigma,\qquad H=G_D^\tau+2\varphi,
 \qquad m_*=\min_F H.
\]
Proposition~\ref{prop:degree} and Theorem~\ref{thm:boundary} give $|\mu|=1$: asymptotically all $n$ pole degrees are used and no mass is lost at the boundary. Hence $|\tau|=|\mu|+|\sigma|=1+|\sigma|<2$, and Lemma~\ref{lem:green-gap} gives $m_*<2\ell_\beta$. Let $\Phi_{\mathrm{sec}}$ map the finite rational coordinate back to the scaled sector variable.  Use $\rho_n=e^{-n\Phi_{\mathrm{sec}}}$, so $\Rea\Phi_{\mathrm{sec}}=2\varphi$; any fixed analytic nonvanishing density factor on a neighborhood of $F$ is also allowed.  There exist $p_n,q_n$ and one table $\omega_{N_n}$, with its zeros on $E$, such that
\begin{equation}\label{eq:sync-degrees}
 \begin{split}
 \deg p_n/n&\longrightarrow|\sigma|<1,
       \qquad \deg q_n\le k_n=o(n),\\
 N_n&=d_n+\deg p_n+k_n+2
       \ge d_n+\deg(p_nq_n)+2,\\
 n^{-1}\mathsf Z(\omega_{N_n})&\weakto\tau^E.
 \end{split}
\end{equation}
With one common scalar normalization, this same table and test satisfy
\begin{equation}\label{eq:sync-upper-multiplier}
 \limsup n^{-1}\log
 \left\|\frac{Q_np_nq_n}{\omega_{N_n}}\right\|_K\le0
       \qquad(K\Subset D),
\end{equation}
and
\begin{equation}\label{eq:sync-lower}
 \left|\frac1{2\pi i}\int_F
 \frac{Q_n(t)p_n(t)q_n(t)\rho_n(t)\chi(t)}
      {\omega_{N_n}(t)}\,dt\right|^{1/n}\longrightarrow e^{-m_*}.
\end{equation}
\end{theorem}

\begin{proof}
Choose the mirror polynomial $p_n$ and the base table $\omega_n^0$ from Proposition~\ref{prop:base-table}, fixing their normalization at that stage. Put
\[
 J_n(z)=\frac{Q_n(z)p_n(z)\rho_n(z)\chi(z)}{\omega_n^0(z)}.
\]
The logarithmic exponent of this same $J_n$ is the sum of the Green contribution from the base table and the external-field contribution from the varying density:
\[
 \frac1n\log|J_n|
 =\frac1n\log\left|\frac{Q_np_n}{\omega_n^0}\right|
   +\frac1n\log|\rho_n|+o(1)
 \longrightarrow -G_D^\tau-2\varphi=-H
\]
in local planar $L^1$; Lemma~\ref{lem:arc-log-means} gives the corresponding statement on small regular arc segments.  Thus the finite table realizes precisely the potential $H$ selected earlier, including its scalar normalization.  For a compact $K\subset F$, choose near-maximizing points and a convergent subsequence. Equation~\eqref{eq:base-moving-envelope}, continuity of $\varphi$, and boundedness of the fixed density give
\[
 \limsup_n n^{-1}\log\|J_n\|_K\le-\inf_K H.
\]
The exact local mirror zeros give \eqref{eq:paired-factorization}; all remaining factors, including the change-of-variable derivative, are analytic and nonvanishing near $z_0$ after shrinking if needed.

Use the boundary point $a\in\partial E$ already used to append zeros to the base table, and apply Lemma~\ref{lem:rational-phase} to the fixed integrand $J_n$. Define
\begin{equation}\label{eq:sync-pair}
 \omega_{N_n}(z)=\omega_n^0(z)(z-a)^{k_n},
 \qquad s_n(z)=p_n(z)q_n(z).
\end{equation}
By construction,
\begin{equation}\label{eq:common-integrand}
 \frac{Q_np_nq_n\rho_n\chi}{\omega_{N_n}}
       =\vartheta_nJ_n
\end{equation}
is an exact identity: the denominator of $\vartheta_n$ is exactly the factor inserted in the final table. The degree relations \eqref{eq:sync-degrees} follow, and \eqref{eq:sync-upper-multiplier} follows from \eqref{eq:base-table-limits} and the global bound for $\vartheta_n$.

On the local segment $J=\psi(I)$, Lemma~\ref{lem:arc-log-means} and finiteness of $H(z_0)$ give
\[
 \left(\int_I|J_n(\psi(x))\psi'(x)|\,dx\right)^{1/n}
       \longrightarrow e^{-m_*}.
\]
We now rule out cancellation for the \emph{same} integrand in \eqref{eq:common-integrand}, first on the local arc and then against the remainder of the contour.  Orient $J$ by increasing $x$. Equation~\eqref{eq:positive-local-integral} and the local Laplace limit give
\begin{equation}\label{eq:sync-local-lower-explicit}
 \Rea\int_J\vartheta_n J_n\,dz
 \ge(1-o(1))\int_J|J_n\,dz|
 =e^{-nm_*+o(n)}.
\end{equation}
Since $H$ is lower semicontinuous on the compact $F$ and has $z_0$ as its unique global minimizer (Theorem~\ref{thm:green-GL9}), the closed set $F\setminus J^\circ$, where $J^\circ=\psi(\operatorname{int}I)$, satisfies
\[
 \inf_{F\setminus J^\circ}H\ge m_*+\delta
 \qquad\text{for some }\delta>0.
\]
The moving-point upper envelope for $J_n$, the fixed contour length, and $\|\vartheta_n\|_F\le e^{o(n)}$ therefore yield
\[
 \left|\int_{F\setminus J}\vartheta_n J_n\,dz\right|
 \le e^{-n(m_*+\delta)+o(n)}.
\]
The exterior contribution is therefore exponentially smaller than the sign-corrected local contribution. Applying the reverse triangle inequality to these two pieces gives
\[
 \left|\int_F\vartheta_n J_n\,dz\right|
 \ge e^{-nm_*+o(n)}-e^{-n(m_*+\delta)+o(n)}.
\]
Conversely, the upper envelope on all of $F$ and the global subexponential bound for $\vartheta_n$ give the upper root bound $e^{-m_*}$.  The fixed factor $(2\pi i)^{-1}$ does not affect the root limit. Thus \eqref{eq:sync-lower} holds without noncancellation of the uncorrected integral.
\end{proof}

\begin{remark}[Why the synchronization is exact]
The table may depend on $Q_n$; what matters is that the residual identity and phase argument use the same table and integrand for each $n$. Although $\deg s_n<n$ eventually, the Hermite identity itself only requires $N_n\ge d_n+\deg s_n+2$; the table degree need not be bounded by $n$. The strict bound $|\sigma|<1$ is instead needed for $|\mu+\sigma|<2$ and hence the exponent gap.
\end{remark}

\begin{proposition}[The critical residual upper bound]
\label{prop:sync-residual-upper}
For a reduced root-optimal sequence in the sector problem, the functional in \eqref{eq:sync-lower} satisfies
\begin{equation}\label{eq:sync-critical-upper}
 \limsup|\mathcal L_n^{\mathrm{ad}}(s_n)|^{1/n}\le e^{-2\ell_\beta}.
\end{equation}
\end{proposition}

\begin{proof}
Fix the tail truncation with error $O(e^{-\kappa_{\mathrm{tail}}n})$, $\kappa_{\mathrm{tail}}>2\ell_\beta$, and all objects of Theorem~\ref{thm:adaptive-sync} on that contour. The regularized target $\mathcal F_n$ is a constant plus a finite Cauchy transform; \eqref{eq:sync-degrees} therefore permits the exact Hermite identity for \eqref{eq:common-integrand}. The tolerance introduced below may change the separating contour, but neither the table nor the integrand defining this functional.

Let $\mathcal M$ be the fixed M\"obius map from the lens chart of Theorem~\ref{thm:boundary} to the current finite chart. To distinguish the two charts, write $\mathcal F_n^{\mathrm{lens}}$ for the truncated target denoted $F_n$ in Lemma~\ref{lem:regularize}. The same target and approximant are represented in the finite chart by $\mathcal F_n=\mathcal F_n^{\mathrm{lens}}\circ\mathcal M^{-1}$ and $R_n=\mathcal R_n\circ\mathcal M^{-1}$. For fixed $\epsilon>0$, the two-constants estimate and logarithmic-radius selection leading to \eqref{eq:error-contour} give lens-chart curves $\Gamma_{t_n}$, $t_n\in[s_1,s_2]$, with
\[
 \|\mathcal F_n^{\mathrm{lens}}-\mathcal R_n\|_{\Gamma_{t_n}}
       \le\exp\{(-2\ell_\beta+\epsilon)n+o(n)\}.
\]
Choose the outer height $s_2$ small enough that the closed collar avoids both the truncated contour and the pole of $\mathcal M$. The image curves $\widetilde\Gamma_{t_n}=\mathcal M(\Gamma_{t_n})$ then separate $E$ from $F$ and lie in a fixed compact band in $D$. They need not be Green level curves for the new point at infinity. Composition preserves the displayed error bound. Their lengths are uniformly bounded because the old curves range over a compact conformal annulus and $\mathcal M'$ is bounded there. Equation~\eqref{eq:sync-upper-multiplier} bounds the multiplier on this entire image band by $e^{o(n)}$. The Hermite identity on these same curves therefore gives
\[
 \limsup|\mathcal L_n^{\mathrm{ad}}(s_n)|^{1/n}
       \le e^{-2\ell_\beta+\epsilon}.
\]
Crucially, $\epsilon$ changes only the contour used to estimate the already fixed functional $\mathcal L_n^{\mathrm{ad}}$. Letting $\epsilon\downarrow0$ therefore proves the bound for exactly the functional appearing in \eqref{eq:sync-lower}.
\end{proof}

\begin{proof}[Proof of Theorem~\ref{thm:root-optimal-rigidity}]
By degree usage, the total normalized pole mass tends to one. Compactness of the sphere gives subsequential weak limits, and Theorem~\ref{thm:boundary} excludes boundary mass, so each limit is a probability measure in $D$.  Suppose one of them is $\mu\ne\lambda_F$.  Fix the tail truncation and apply Theorem~\ref{thm:green-GL9}, Theorem~\ref{thm:adaptive-sync}, and Proposition~\ref{prop:sync-residual-upper} along this subsequence.  They give, for exactly the same functional,
\[
 e^{-m_*}
 =\lim|\mathcal L_n^{\mathrm{ad}}(s_n)|^{1/n}
 \le e^{-2\ell_\beta},\qquad m_*<2\ell_\beta,
\]
a contradiction.  Thus every subsequential limit on the compactified sphere is $\lambda_F$. The inverse M\"obius map is continuous on the sphere, and Theorem~\ref{thm:boundary} excludes limiting mass at the compactified boundary, including the points corresponding to zero and infinity. Pushing the measures back and undoing the scale $hz/n$ gives \eqref{eq:all-root-optimal-poles} with no lost mass. The root theorem applies to best approximants, proving the last assertion.
\end{proof}

\begin{remark}[Scope and logical dependencies]
Nonseparation, the S-variational conditions, and the target-specific boundary argument are essential; arbitrary separating supports are not covered. The allowance for $o(n)$ exceptional poles is also necessary: \ref{app:exceptional-poles} inserts prescribed poles with tiny nonzero residues while preserving root optimality. Their normalized number can tend to zero arbitrarily slowly. Root optimality alone therefore yields no uniform rate of weak convergence for the pole measures.
\end{remark}

\section{Constructive consequences}\label{sec:consequences}

The weak pole law identifies the limiting distribution, but neither fixes individual pole locations nor constructs approximants with a common denominator. This section gives two complementary existence results: a scalar-to-joint error comparison for finite $\varphi$-families, and a selection of root-optimal families whose poles all escape fixed compact sets. We take $h=1$; dilation restores general $h$. Neither construction is needed for the rigidity proof.

\subsection{One denominator for a finite family of \texorpdfstring{$\varphi$}{phi}-functions}\label{subsec:common-denominator}
For $j\ge1$ define
\[
 \varphi_j(w)=\frac{e^w-\sum_{k=0}^{j-1}w^k/k!}{w^j},
 \qquad \varphi_0(w)=e^w,
\]
with the removable value at $w=0$. These functions occur in exponential integrators \cite{HochbruckOstermann2010}, where shared poles are useful \cite{SchmelzerTrefethen2007,AlMohy2025}. We compare the best joint and scalar errors on the whole infinite sector.

A numerator Taylor correction preserves the denominator; see also \cite[\S5]{Schmelzer2026}. To bound the Taylor correction, however, we need derivative estimates at the vertex, and a sector bound alone gives no full disk there. A positive translation supplies one: for $0<\beta<\pi$ and $a>0$, we have $a+\Sigma_\beta\subset\Sigma_\beta$, and the disk of radius $a\sin(\beta/2)$ centered at $a$ lies inside the sector. After translation, the error is bounded on $\Sigma_\beta$ and on a full disk about zero.

\begin{proposition}[Translation and Taylor correction]\label{prop:translation-jets}
Let $0<\beta<\pi$, $r=P/Q\in\RR_n(\Sigma_\beta)$ be reduced, and
\[
 e=\|e^{-z}-r(z)\|_{\Sigma_\beta}.
\]
Fix $a>0$ and $0<R<a\sin(\beta/2)$. Set
\begin{equation}\label{eq:translation-jets}
 s(z)=e^a r(z+a),\qquad s_0=s,\qquad
 s_j(z)=\frac{s(z)-T_{j-1}s(z)}{(-z)^j}\quad(j\ge1),
\end{equation}
where $T_{j-1}s$ is the Taylor polynomial of $s$ at zero through degree $j-1$. Every $s_j$ has a representation with denominator $Q(z+a)$. For $j\ge1$, its numerator degree is at most $n-1$. Moreover,
\begin{align}
 \|e^{-z}-s_0\|_{\Sigma_\beta}&\le e^a e,\label{eq:translation-exp}\\
 \|\varphi_j(-z)-s_j\|_{\Sigma_\beta}
 &\le (3\cdot2^j-2)e^aR^{-j}e\qquad(j\ge1).
 \label{eq:translation-phi}
\end{align}
All poles of these functions have distance at least $a\sin(\beta/2)$ from $\Sigma_\beta$ in the finite plane.
\end{proposition}

\begin{proof}
Write $\Delta(z)=e^{-z}-s(z)$ and $M=e^ae$. The translation observation gives $|\Delta|\le M$ on $\Sigma_\beta$ and on the closed disk $|z|\le R$. If $\Delta(z)=\sum_{k\ge0}d_kz^k$ near zero, Cauchy's estimate gives $|d_k|\le MR^{-k}$. Since
\[
 \varphi_j(-z)-s_j(z)=\frac{\Delta(z)-T_{j-1}\Delta(z)}{(-z)^j},
\]
the Taylor tail on $|z|\le R/2$ is bounded by $2MR^{-j}$. For $z\in\Sigma_\beta$ with $|z|\ge R/2$, the triangle inequality instead gives
\[
 \frac{|\Delta(z)|}{|z|^j}
 +\sum_{k=0}^{j-1}|d_k|\,|z|^{k-j}
 \le MR^{-j}\left(2^j+\sum_{k=0}^{j-1}2^{j-k}\right)
 =(3\cdot2^j-2)MR^{-j}.
\]
This proves the error estimates, including the removable value at zero.

To check the rational type, let $m=\deg Q\le n$. The numerator of
\[
 s-T_{j-1}s
 =\frac{e^aP(z+a)-Q(z+a)T_{j-1}s(z)}{Q(z+a)}
\]
has a zero of order at least $j$ at zero, and degree at most $m+j-1$. Dividing it by $(-z)^j$ leaves a polynomial of degree at most $m-1$; no new pole is created at the origin. The zero function covers the case $m=0$.

Finally put $\theta=\beta/2$. The sector is the intersection of the two half-planes $L_\pm(z)\ge0$, where
\[
 L_\pm(z)=\sin\theta\,\Rea z\ \pm\cos\theta\,\Ima z.
\]
These are signed distances because their gradients have length one. For every original pole $p\notin\Sigma_\beta$, at least one $L_\pm(p)$ is negative, and then $L_\pm(p-a)<-a\sin\theta$. The poles of the translated functions are among these $p-a$, proving the separation claim.
\end{proof}

For a fixed integer $L\ge0$, let $\mathcal E^{\rm com}_{n,L}(\beta)$ be the infimum of
\[
 \max_{0\le j\le L}\left\|\varphi_j(-z)-\frac{P_j(z)}{Q(z)}\right\|_{\Sigma_\beta}
\]
over a single nonzero denominator $Q$, with no finite zeros on $\Sigma_\beta$, and polynomials satisfying $\deg Q\le n$ and $\deg P_j\le\deg Q$. Cancellations in individual components are allowed.

\begin{corollary}[Degree-independent common-denominator comparison]\label{cor:common-denominator}
For every fixed $0<\beta<\pi$ and $L\ge0$ there is $C_{\beta,L}<\infty$, independent of $n$, such that
\begin{equation}\label{eq:common-comparison}
 \eps_n(\beta;1)\le\mathcal E^{\rm com}_{n,L}(\beta)
                 \le C_{\beta,L}\eps_n(\beta;1).
\end{equation}
Consequently $(\mathcal E^{\rm com}_{n,L}(\beta))^{1/n}\to\rho_\beta$.
If a sequence of such common-denominator families has joint error with root limit $\rho_\beta$, its common denominator, after factors common to all numerators are removed, has scaled zero measure converging to $\nu_\beta$.
\end{corollary}

\begin{proof}
The lower bound follows by retaining the exponential component $j=0$. For the upper bound, use a scalar best approximant in Proposition~\ref{prop:translation-jets}, with fixed $a$ and $R$. One may take
\[
 C_{\beta,L}=e^a\max\!\left\{1,\max_{1\le j\le L}(3\cdot2^j-2)R^{-j}\right\},
\]
omitting the inner maximum for $L=0$. Taking roots gives the rate.
For the pole statement, consider the component $j=0$, which approximates the exponential with root-optimal error. After reduction, Proposition~\ref{prop:degree} gives it $n-o(n)$ poles, and Theorem~\ref{thm:root-optimal-rigidity} gives its limiting pole measure. Since this reduced denominator divides the common denominator of degree at most $n$, at most $o(n)$ additional zeros remain. Their normalized mass vanishes, so the common denominator has the same weak limit.
\end{proof}

\subsection{Selecting root-optimal sequences without finite exceptional poles}
Weak convergence permits exceptional poles, so it cannot require every pole to leave a fixed compact set. The translation argument nevertheless lets us \emph{select} a root-optimal sequence with that stronger property.

\begin{corollary}[Escaping poles and local uniform approximation]\label{cor:escaping-selection}
Fix $0<\beta<\pi$ and $L\ge0$. There exist common-denominator families $s_{j,n}$, $0\le j\le L$, of denominator degree at most $n$, whose joint error on $\Sigma_\beta$ has root limit $\rho_\beta$, whose poles all escape every fixed compact subset of $\C$, and such that, for every compact $K\subset\C$,
\begin{equation}\label{eq:selected-overconvergence}
 \limsup_{n\to\infty}
 \max_{0\le j\le L}\|\varphi_j(-z)-s_{j,n}\|_K^{1/n}
 \le\rho_\beta.
\end{equation}
Their common-denominator scaled pole measures still converge to $\nu_\beta$.
\end{corollary}

\begin{proof}
Start with scalar best approximants and apply Proposition~\ref{prop:translation-jets} with $a=a_n\to\infty$, $a_n=o(n)$, and $R_n=a_n\sin(\beta/2)/2$. The constants in the error bounds have logarithm $o(n)$, while the exponential component supplies the lower bound $\eps_n(\beta;1)$. Thus the joint root limit is $\rho_\beta$. All poles have distance at least $a_n\sin(\beta/2)$ from the sector, so, because $0\in\Sigma_\beta$, their moduli tend to infinity uniformly. Every fixed $K$ eventually lies in $|z|\le R_n/2$. The disk estimates in the proof of Proposition~\ref{prop:translation-jets} then give \eqref{eq:selected-overconvergence}. The common-denominator pole law follows from Corollary~\ref{cor:common-denominator}. In particular, the sublinear translation does not alter the limiting distribution on the degree scale.
\end{proof}

The positive opening supplies the disk used for Cauchy estimates. For example, $a=1$ and $R=\sin(\beta/2)/2$ give $C_{\beta,L}=O_L(\beta^{-L})$ for fixed $L\ge1$ as $\beta\downarrow0$. This is a bound for this construction, not a necessary blow-up of the best comparison constant. At $\beta=0$ the disk disappears, so the argument does not settle the half-line questions of \cite{StahlSchmelzer2009,Schmelzer2026}. Pole escape is asserted for the selected families, not for exact best approximants.

\subsection{Why a single repeated pole cannot attain the free-pole rate}
The next comparison illustrates how rigidity can yield a strict inequality between convergence factors without estimating their difference directly. The repeated-pole construction of \cite[Theorem~4.5 and Proposition~4.6]{XueZhangRepeated2026} attains
\[
 U_\beta=\tan\!\left(\frac{\pi}{4(2-\beta/\pi)}\right)
\]
when its pole scale is optimized within the noncoalescing branch analyzed there. It is a competitor for the free-pole problem, so $\rho_\beta\le U_\beta$. Rigidity makes the inequality strict.

\begin{corollary}[Strict gain from free poles]\label{cor:strict-free-rate}
For every $0\le\beta<\pi$, $\rho_\beta<U_\beta$ for the attained repeated-pole factor just defined.
\end{corollary}

\begin{proof}
Put $\gamma=2-\beta/\pi$ and $\alpha_\beta^*=1/[\gamma\sin(\pi/(2\gamma))]$.  The repeated-pole result in \cite[Theorem~4.5 and Proposition~4.6]{XueZhangRepeated2026} supplies a sequence with sole possible pole $-\alpha_\beta^*n/h$ and root error $U_\beta$.  If $\rho_\beta=U_\beta$, Proposition~\ref{prop:degree} makes its reduced pole mass tend to one.  Its scaled pole measures then converge to $\delta_{-\alpha_\beta^*}$, whereas Theorem~\ref{thm:root-optimal-rigidity} gives the non-atomic arc measure $\nu_\beta$, a contradiction.
\end{proof}

\section{Numerical tests of the rate and pole geometry}\label{sec:num}

The experiments compare numerical candidates with freely fitted poles against the classical root factor and the predicted pole measure. Table~\ref{tab:num} summarizes both comparisons; Supplementary Material, Section S2, provides further scalar checks, angle-dependent results, and repeated-pole comparisons.

We set $h=1$ and fit general complex type-$(n,n)$ candidates by AAA followed by Lawson iteration \cite{NakatsukasaSeteTrefethen2018,NakatsukasaTrefethen2020}, using boundary-ray samples in the coordinate $w=(n-z)/(n+z)$. Neither the S-curve nor predicted poles enter the fit, and no pole--zero pairs are removed. Independent denser boundary grids and local peak searches estimate the error. The precision-sensitive case $\beta=120^\circ$, $n=24$ is evaluated at the refined peaks in 80-digit arithmetic, with a 100-digit check of the largest peak. Each candidate has $n$ extracted finite poles and no detected poles in $\Sigma_\beta$. Refined residual searches check the estimated errors, while denominator winding counts provide an independent numerical check for sector poles. Supplementary Material, Section S2, records the protocol; the reproduction archive contains the scripts and stored candidates. This finite-step AAA--Lawson computation has no convergence guarantee to a local or global minimizer on the continuous sector.

For the geometric comparison, let $z_{j,n}=p_{j,n}/n$ and $\Gamma_\beta=\operatorname{supp}\nu_\beta$. For each $z_{j,n}$, choose a nearest point on the polygonal curve joining $65537$ equally spaced equilibrium-mass nodes on $\Gamma_\beta$.  Write $d_j$ for the Euclidean distance to that curve, $L_\beta$ for the Euclidean distance between its endpoints, and $q_j\in[0,1]$ for the interpolated cumulative equilibrium mass measured from the lower endpoint.  Let $q_{(j)}$ denote the sorted mass coordinates.  We report the $90$th percentile $D_{90}$ of $d_j/L_\beta$ and the midpoint-quantile discrepancy
\[
 D_{\mathrm{q}}=\frac1n\sum_{j=1}^n\left|q_{(j)}-\frac{j-1/2}{n}\right|.
\]
$D_{\mathrm q}$ measures discrepancies in the projected mass coordinates, while $D_{90}$ disregards the largest $10\%$ of distances; neither diagnostic certifies weak convergence. Supplementary Material, Section~S2, also reports an all-pole mean spherical-chordal deviation, which does not discard a fixed fraction of the poles. Doubling the curve resolution changes the original two diagnostics by less than $4\times10^{-9}$.

\begin{table}[H]
\caption{Unrestricted AAA--Lawson candidates.  Here $E_n$ is the estimated boundary error; $E_n^{1/n}$ is computed before rounding.  Neither sampled errors nor the geometric diagnostics certify finite-degree minimax optimality.}\label{tab:num}
\centering
\small
\begin{tabular}{c r c c c c c}
\toprule
$\beta$ & $n$ & $E_n$ & $E_n^{1/n}$ & $\rho_\beta$ & $D_{90}$ & $D_{\mathrm{q}}$\\
\midrule
$90^\circ$ & 16 & $2.26\times10^{-11}$ & 0.216112 & 0.226059 & 0.006306 & 0.004236\\
$120^\circ$ & 24 & $8.83\times10^{-13}$ & 0.314594 & 0.321876 & 0.003506 & 0.002484\\
$150^\circ$ & 32 & $3.80\times10^{-10}$ & 0.507709 & 0.512291 & 0.001556 & 0.001396\\
$170^\circ$ & 72 & $1.23\times10^{-8}$ & 0.776498 & 0.775255 & 0.000096 & 0.000316\\
$175^\circ$ & 120 & $1.78\times10^{-7}$ & 0.878525 & 0.876345 & 0.000149 & 0.001781\\
$179^\circ$ & 240 & $2.31\times10^{-3}$ & 0.975024 & 0.973161 & 0.000673 & 0.000575\\
\bottomrule
\end{tabular}

\end{table}

\begin{figure}[H]
\centering
\begin{minipage}[t]{.49\textwidth}
\centering
\includegraphics[width=\linewidth]{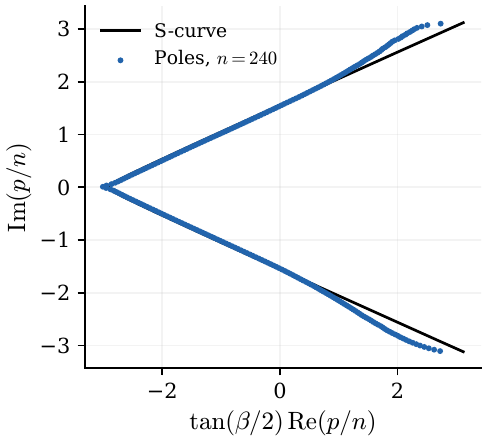}\\[-1mm]
\small (a) Scaled poles and theoretical S-curve.
\end{minipage}\hfill
\begin{minipage}[t]{.49\textwidth}
\centering
\includegraphics[width=\linewidth]{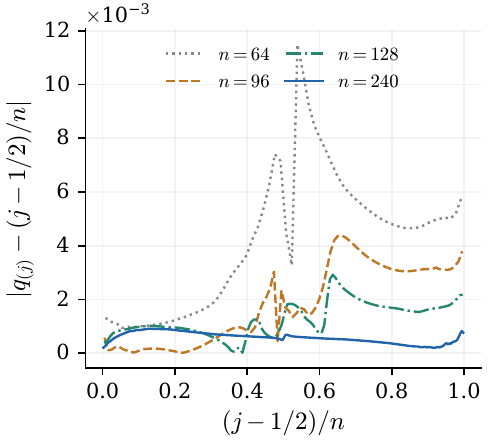}\\[-1mm]
\small (b) Projected mass-coordinate errors.
\end{minipage}
\caption{Pole tests for $\beta=179^\circ$.  Panel (a) compares the $n=240$ poles with the predicted support; the horizontal coordinate is magnified by $\tan(\beta/2)$.  Panel (b) compares the projected mass coordinates for four degrees.}\label{fig:179numerics}
\end{figure}

We compare $E_n^{1/n}$ with the root factor in Theorem~\ref{thm:root}: for $\beta=170^\circ,175^\circ,179^\circ$, their relative differences are below $0.25\%$.  Since both factors are close to one, this comparison compresses differences in the logarithmic decay rate.  At $\beta=179^\circ$ and $n=240$, for example,
\[
 -\frac{\log E_n}{n}\approx0.0252928,
 \qquad -\log\rho_\beta\approx0.0272060.
\]
The estimated finite-degree decay exponent is therefore about $7.03\%$ below the limiting exponent, equivalently $E_n/\rho_\beta^n\approx1.583$. The root limit does not determine whether $E_n^{1/n}$ lies above or below $\rho_\beta$. Nor are the geometric diagnostics expected to be monotone across rows with different $\beta$ and $n$. These tests support the predicted rate and geometry, without certifying root optimality, continuum suprema, or finite-degree minimax optimality.

For $\beta=179^\circ$ and $n=64,96,128,240$, the error estimates are $0.24311$, $0.09763$, $0.04239$, and $0.002311$, respectively.  The corresponding $D_{\mathrm{q}}$ values are $0.004174$, $0.001808$, $0.001265$, and $0.000575$.  The long-interval factor $(E_{240}/E_{64})^{1/176}=0.973892$ is close to $\rho_{179^\circ}=0.973161$, although this secondary diagnostic depends on the chosen degrees.

Certifying finite-degree minimax errors requires further work. The convex-dual/Ruttan framework of \cite{ZhangYangYangZhang2025} offers one possible route, separate from the continuous-sector asymptotic proof.

\section{Conclusions}\label{sec:conclusion}

A real Ferrers zero determines the classical sector root factor and, after mass normalization, the shape and scale of the equilibrium arc. Theorem~\ref{thm:root-optimal-rigidity} then establishes that every reduced root-optimal sequence has this same scaled pole law. The new pole-rigidity range is that of positive-opening sectors; the corresponding half-line law was previously announced in \cite{StahlSchmelzer2009}.

The proof first excludes boundary mass by pole stripping. For an incorrect interior limit, Green reflection and a normalized Hermite test give incompatible exponential bounds on one integral. Separately, translation and Taylor correction give finite $\varphi$-families with a common denominator at the scalar root rate, including selected families whose poles all escape fixed compact sets.

Root optimality alone permits exceptional-pole fractions that tend to zero arbitrarily slowly (\ref{app:exceptional-poles}). Thus it yields no uniform quantitative rate for convergence of the pole measures without further assumptions. Sharp error prefactors, the wide-sector limit of the appropriately rescaled equilibrium measure as $\beta\uparrow\pi$, and stable numerical realization of the common-denominator construction remain open directions.

\appendix
\renewcommand{\thetheorem}{\Alph{section}.\arabic{theorem}}

\section{Analytic identities for the sector equilibrium}\label{app:equilibrium}
Two identities underpin Section~\ref{sec:Scurve}: a first-order equation for the Stieltjes differential, and a half-period formula that determines the mass and equilibrium level. Their statements appear in that section; the calculations follow here.

\subsection{Mellin normalization and the differential equation}\label{app:mellin}
\begin{proof}[Proof of Lemma~\ref{lem:mellin-differential}]
The real-integral normalization can be checked without choosing a trajectory. For $0<s<1$ and $-1<c<1$, the substitution $t=x/(1-x)$ and the binomial expansion, followed by beta integration, give
\begin{equation}\label{eq:mellin-audit}
 \int_0^\infty\frac{t^{s-1}}{\sqrt{t^2+2ct+1}}\,dt
 =\frac{\pi}{\sin(\pi s)}\,
 {}_2F_1(s,1-s;1;(1-c)/2)
 =\frac{\pi}{\sin(\pi s)}\mathsf P_{s-1}(c).
\end{equation}
After expanding, the $k$th term contains the Euler beta integral $\mathrm B(s+k,1-s+k)$. The relation $(2k)!=4^kk!(1/2)_k$ simplifies the binomial coefficient; here $(a)_k=a(a+1)\cdots(a+k-1)$ denotes the rising factorial. Continuing to $s=1+\eta$ by subtracting the leading term at infinity yields
\begin{equation}\label{eq:regularized-mellin}
 J_0:=\int_0^\infty
 \left(\frac{t^\eta}{\sqrt{t^2+2\Lambda ct+\Lambda^2}}-t^{\eta-1}\right)dt
 =-\frac{\pi\Lambda^\eta}{\sin(\pi\eta)}\mathsf P_\eta(c)=0.
\end{equation}
For clarity, the continuation is the sum of the integral over $(0,1)$, the subtracted integral over $(1,\infty)$, and $1/(1-s)$ before rescaling. Both subtracted integrals converge in $0<\Rea s<2$, away from $s=1$.

On $w=-t$, $t>0$, the physical branch has $R(-t)=-\sqrt{t^2+2\Lambda ct+\Lambda^2}$.  The jump of $S$ is $-2\pi i t^\eta/\sqrt{t^2+2\Lambda ct+\Lambda^2}$, which cancels the jump of $\gamma w^\eta/2$ in $B$.  At zero,
\[
 S(w)=S(0)-\frac{\pi}{\Lambda\sin(\pi\eta)}w^\eta+O(w),
 \qquad B(w)=\kappa\Lambda S(0)+O(w),
\]
so the remaining isolated singularity is removable.  Splitting the Stieltjes integral into $t<|w|/2$, $|w|/2<t<2|w|$, and $t>2|w|$, and subtracting the local density in the middle interval, gives $B(w)=O(|w|^\eta)$ at infinity, including its continued values across the negative axis.  Since $\eta<1$ and $B$ is analytic outside a compact set, the Laurent expansion at infinity has no positive powers: $B$ has a removable singularity there. Along $w=x>0$,
\[
 S(x)=\frac{\pi}{\sin(\pi\eta)}x^{\eta-1}+\frac{J_0+o(1)}x.
\]
Thus $B(\infty)=-\kappa J_0=0$, proving $B=O(w^{-1})$.  On the endpoint cut, $R$ changes sign while $S$ is analytic there, so
\[
 B_++B_-=\gamma w^\eta.
\]
Moreover
\[
 B-\frac\gamma2w^\eta=-\kappa RS
\]
vanishes like a square root at $A_\pm$, because $S$ is finite and nonzero there.  This verifies the jump and soft-edge conditions directly.

Put $\mathcal F(w)=B(w)-\gamma w^\eta/2$.  The particular differential combination below is chosen because it annihilates the prescribed $w^\eta$ jump.  Indeed,
\[
 R(w)\{w\mathcal F'(w)-\eta\mathcal F(w)\}
\]
has no jump on the negative axis, while on the endpoint cut the sign changes of $R$ and $\mathcal F$ cancel.  Its endpoint and zero singularities are removable and it is bounded at infinity.  Hence it is a single-valued bounded meromorphic function on the sphere and therefore a constant by Liouville's theorem.  Its value at zero is
\[
 C=\kappa\eta\Lambda^2 S(0)
  =\frac{\gamma\eta}{2}\Lambda^{\eta+1}\mathsf P_{\eta-1}(c)>0,
\]
by \eqref{eq:mellin-audit}.  The Legendre recurrence at $\mathsf P_\eta(c)=0$ gives the second expression in \eqref{eq:audited-ODE}; setting $w=u^2$ proves the differential identity. Finally, $\Ima S(w)$ has the opposite sign to $\Ima w$, and $S(w)>0$ for $w>0$. Thus it has no zero in $\C\setminus(-\infty,0]$, as claimed.
\end{proof}

\subsection{The branch-explicit half-period}\label{app:period}
\begin{proof}[Proof of Lemma~\ref{lem:half-period-audit}]
Set $u=\sqrt\Lambda\tan v$ and $k=\cos(\theta/2)=\sqrt{1-m}$. The differential becomes
\[
 \frac{du}{\widehat r(u)}=\frac1{2\sqrt\Lambda}
       \frac{d(2v)}{\sqrt{1-k^2\sin^2(2v)}}.
\]
At $u=a_+$, $2v=\pi/2+ib$, where $\cosh b=1/k$.  Continue first along the real segment to $\pi/2$ and then vertically.  The real part is $\Kell(k^2)$; the imaginary part is
\[
 \int_0^b\frac{dt}{\sqrt{1-k^2\cosh^2t}}
 =\int_0^{\pi/2}\frac{d\psi}{\sqrt{k^2+m\sin^2\psi}}=\Kell(m),
\]
using $\sinh t=(\sqrt m/k)\sin\psi$.  This proves the half-period with its signs, using the elliptic \emph{parameter} convention. Integrate \eqref{eq:audited-ODE} and use $\widehat F(a_+)=T(a_+)=0$ to obtain \eqref{eq:T-audit}.
\end{proof}

\section{Proof of the theta--Ferrers correspondence}\label{app:theta}

GR's sector formula \cite[\S2.5, p.~329, Eq.~(27)]{GoncharRakhmanov1989} is
\[
 \int_0^{1/2}\left(\frac{\vartheta_0(t)}{\vartheta_3(t)}\right)^{2\alpha_{\rm GR}}dt=0,
 \qquad \alpha_{\rm GR}=1-\frac{\theta_{\rm GR}}\pi,
 \qquad v_{\theta_{\rm GR}}=-\mathfrak h^2.
\]
The dictionary to our notation is
\[
 \theta_{\rm GR}=\beta/2,\quad \alpha_{\rm GR}=\eta+\tfrac12,\quad
 \mathfrak h=i\sqrt q,\quad
 \vartheta_0(t)=\theta_4(\pi t,\mathfrak h),\quad
 \vartheta_3(t)=\theta_3(\pi t,\mathfrak h).
\]
Thus $v_{\theta_{\rm GR}}=q$: the imaginary nome is not the positive root factor. We now prove the identity in Proposition~\ref{prop:theta-ferrers}, with the standard theta and Jacobi conventions of \cite[\S\S20.7, 22.2]{DLMF}.

\begin{proof}[Proof of Proposition~\ref{prop:theta-ferrers}]
Fix $q\in(0,1)$ and put
\[
 a_0=\theta_3(0,q^2),\quad b_0=\theta_2(0,q^2),\quad
 r=b_0/a_0\in(0,1),\quad T=a_0^2+b_0^2.
\]
Separating even and odd theta indices gives
\[
 A:=\theta_4(0,i\sqrt q)=a_0-ib_0,\qquad
 B:=\theta_3(0,i\sqrt q)=a_0+ib_0.
\]
The duplication identities yield
\begin{equation}\label{eq:theta-duplication}
 AB=T=\theta_3(0,q)^2=\frac{2\Kell(m)}\pi,
 \qquad \theta_2(0,q)^2=2a_0b_0,
\end{equation}
where $m=m(q)$. Set $\vartheta=4\arctan r\in(0,\pi)$. Then
\[
 A/B=e^{-i\vartheta/2},\qquad
 \cos\vartheta=1-\frac{8r^2}{(1+r^2)^2}=1-2m.
\]

The same index separation at a variable argument gives
\[
 v(t):=\frac{\theta_4(\pi t,i\sqrt q)}{\theta_3(\pi t,i\sqrt q)}
 =\frac{1-ix(t)}{1+ix(t)},\qquad
 x(t)=\frac{\theta_2(2\pi t,q^2)}{\theta_3(2\pi t,q^2)}.
\]
In this paragraph only, let $k=r^2$ denote a Jacobi \emph{modulus}, and set $K_2=\Kell(k^2)=\pi a_0^2/2$. The theta representation gives
\[
 x(t)=r\,\operatorname{cd}(4K_2t,k),\qquad
 (x'(t))^2=16K_2^2(r^2-x(t)^2)(1-r^2x(t)^2).
\]
Here $\operatorname{cd}=\operatorname{cn}/\operatorname{dn}$. Its derivative is $-(1-k^2)\operatorname{sn}/\operatorname{dn}^2$, so $x(t)$ decreases from $r$ to $-r$ on $[0,1/2]$. Consequently $v(t)=e^{is(t)/2}$ with $s(t)=-4\arctan x(t)$ increasing from $-\vartheta$ to $\vartheta$. This proves the branch assertion.

Direct algebra now gives
\[
 \cos s-\cos\vartheta
 =\frac{8(r^2-x^2)(1-r^2x^2)}{(1+r^2)^2(1+x^2)^2},
\]
and, taking the positive square root along the increasing path,
\begin{equation}\label{eq:theta-measure}
 dt=\frac{ds}{2\pi T\sqrt{2(\cos s-\cos\vartheta)}}.
\end{equation}
The endpoint singularities are integrable. Substitution and symmetry turn the integral in \eqref{eq:theta-ferrers} into
\[
 \frac{1}{\pi T\sqrt2}\int_0^\vartheta
 \frac{\cos((\eta+1/2)s)}{\sqrt{\cos s-\cos\vartheta}}\,ds.
\]
The Mehler--Dirichlet representation \cite[Eq.~14.12.1]{DLMF} identifies this as $\mathsf P_\eta(\cos\vartheta)/(2T)$. Equation~\eqref{eq:theta-duplication} proves the claimed prefactor.
\end{proof}

At the half-line endpoint $\beta=0$, retain $q,m,T,A,\vartheta$ from the preceding proof with $q=v=\rho_0$. The physical normalization agrees with GR's Lemma~5 as well. Here $\eta=1/2$, $\gamma=2$, and the Legendre recurrence at the root gives
\[
 -\frac32\mathsf P_{3/2}(c)=\frac12\mathsf P_{-1/2}(c)
 =\frac{\Kell(m)}\pi.
\]
Thus \eqref{eq:lambda-scale} reduces to $\Lambda=\pi^2/\Kell(m)^2=4/T^2$. The lower endpoint in the original sector coordinate is
\[
 -\Lambda e^{i\vartheta}=-4A^{-4}
 =-4\theta_4(0,i\sqrt q)^{-4},
\]
which is GR's endpoint formula \cite[Lemma~5, p.~324]{GoncharRakhmanov1989}. The scale as well as the root factor is therefore calibrated against the classical construction.

For comparison, Aptekarev's half-line result \cite{Aptekarev2002} is the stronger statement
\[
 \mathfrak e_n(e^{-x};[0,\infty))=2\sqrt v\,v^n(1+o(1)),\qquad v=\rho_0.
\]
Equivalently, the $2n$th-root limit of this best error is $\sqrt v$, whereas our $n$th-root factor is $v$. Strong asymptotics for a varying exponential on a bounded real interval are treated in \cite{MagnusMeinguet2024}. Neither prefactor statement is being extended by the root-rate identification above.

\section{Auxiliary potential-theoretic estimates}\label{app:potential}
These estimates serve different steps of the rigidity proof. Their hypotheses are stated where used, except for the uniqueness lemma below.

\subsection{Choosing a contour with small total pole potential}\label{app:contour}
\begin{proof}[Proof of Lemma~\ref{lem:cap-select}]
Let $\Psi$ map the exterior of $K$ conformally onto $|\zeta|>1$, fixing infinity.  For a finite pole with height $a=\log|\Psi(p)|>0$, the explicit disk Green kernel gives
\begin{align*}
 &\sup_{z\in\Gamma_s}g_K(z,p)
 =\log\frac{e^{s+a}-1}{|e^s-e^a|}\\
 =&\min(s,a)+\log\frac{1-e^{-(s+a)}}{1-e^{-|s-a|}}\le s+\log\frac{s+a}{|s-a|}.
\end{align*}
The maximum is attained at the angle of $\Psi(p)$; the last inequality uses that $(1-e^{-t})/t$ decreases for $t>0$.  The singularity at $s=a$ is integrable.  Scaling by $a$, splitting the integral at one, and expanding $\log((1+x)/(1-x))$ on $(0,1)$ yield
\begin{equation}\label{eq:log-radius-integral}
 \int_0^\infty\log\frac{s+a}{|s-a|}\,\frac{ds}{s}
 =4\sum_{k=0}^\infty\frac1{(2k+1)^2}=\frac{\pi^2}{2}.
\end{equation}
For a pole at infinity the Green kernel is exactly $s$.  Subadditivity of the maximum and integration with respect to $ds/s$ therefore give
\[
 \int_{s_1}^{s_2}\max_{\Gamma_s}\sum_jg_K(z,p_j)\,\frac{ds}{s}
 \le N(s_2+\pi^2/2).
\]
The averaging principle proves \eqref{eq:log-radius-selection}; its finite right-hand side excludes any pole on the selected curve.  Neither condenser-capacity normalization nor residue estimates are needed.
\end{proof}

\subsection{A separating perturbation with a unique minimizer}\label{app:generic}
\begin{lemma}[Generic uniqueness from separating perturbations]
\label{lem:generic-unique}
Let $X$ be a compact metric space, let $\Phi:X\to(-\infty,+\infty]$ be lower semicontinuous with finite minimum, and let $f_j\in C(X;\R)$ be uniformly bounded and separate points of $X$.  Suppose $a_j>0$ and $\sum_j a_j<\infty$.  Then there are parameters $\alpha_j\in[0,a_j]$ for which
\[
 x\longmapsto \Phi(x)+\sum_{j\ge1}\alpha_j f_j(x)
\]
has a unique minimizer.
\end{lemma}

\begin{proof}
The series converges uniformly, so every perturbed function is lower semicontinuous and attains its minimum.  Give each interval $[0,a_j]$ normalized Lebesgue measure and their countable product the product probability measure. For parameter vectors $\alpha,\widetilde\alpha$, their minimum values satisfy
\[
 |M(\alpha)-M(\widetilde\alpha)|
 \le \bigl(\sup_j\|f_j\|_X\bigr)
           \sum_j|\alpha_j-\widetilde\alpha_j|.
\]
The metric given by this sum induces the product topology because $\sum_j a_j<\infty$.  Hence $M$ is continuous.  Its one-sided partial derivatives are limits of measurable difference quotients, so their nonagreement sets are measurable; parameter endpoints form null sets. Thus the use of Fubini below is justified. Fix all coordinates except $\alpha_j$ and write
\[
 M_j(s)=\min_{x\in X}\{\Phi_j(x)+s f_j(x)\},\qquad 0\le s\le a_j,
\]
where $\Phi_j$ contains the other parameters.  This is the infimum of affine functions of $s$ with uniformly bounded slopes, hence is finite, concave, and Lipschitz.  It is differentiable for almost every $s$.

At a differentiability point, any two minimizers $x,y$ must satisfy $f_j(x)=f_j(y)$.  Indeed, evaluating the minimum at $x$ after a positive or negative increment gives the two one-sided slope bounds
\[
 M'_{j,+}(s)\le f_j(x)\le M'_{j,-}(s),
\]
and the same bounds hold for $y$; differentiability makes both endpoints equal. By Fubini's theorem, for each fixed $j$ this differentiability property holds for product-almost every parameter vector.  Intersecting these full-measure sets over the countable set of indices gives one vector for which every pair of minimizers has the same value of every $f_j$.  Since the $f_j$ separate points, the minimizer is unique.
\end{proof}

\subsection{The spherical balayage constant}\label{app:balayage-constant}
\begin{proof}[Proof of Lemma~\ref{lem:spherical-balayage}]
For the ordinary logarithmic potential $V^\alpha$ defined in Proposition~\ref{prop:balayage}, the point-charge identity is
\[
 V^{\delta_a}(z)-V^{\omega_a}(z)
       =g_D(z,a)-g_D(a,\infty).
\]
It follows from Green representation, or by matching its logarithmic singularity and zero boundary values after adding $g_D(a,\infty)$. Replacing the logarithmic kernel by the chordal kernel adds $L(a)$ to the point charge and $\int L\,d\omega_a$ to its sweep; the common $L(z)$ terms cancel.  This proves the point-charge case of \eqref{eq:spherical-balayage}. Regularity gives weak continuity of $a\mapsto\omega_a$ up to the boundary. Also
\[
 g_D(a,\infty)=\log|a|-\log\capc(E)+o(1)
       \quad(a\to\infty).
\]
These observations prove the stated continuous extensions and boundedness. Integration proves \eqref{eq:spherical-balayage}; the swept potential is finite on every compact subset of $D$, the chordal kernel is nonnegative, and the constant is bounded.  Thus no subtraction of two undefined infinite logarithmic moments is involved.
\end{proof}

\subsection{Logarithmic means on a regular arc}\label{app:arc-means}
\begin{proof}[Proof of Lemma~\ref{lem:arc-log-means}]
For a sufficiently short embedded arc, the set of parameters with $[\psi(x),a]<s$ has length at most $Cs$, uniformly in the source $a\in\Chat$.  Layer integration therefore bounds the integral of the logarithmic singularity on this set by $Cs(1+|\log s|)$.  The source $a=\infty$ is harmless because the chordal kernel is continuous there on the finite arc.  Truncated kernels are continuous on the compact sphere, and this uniform tail bound proves continuity of the full kernel as an $L^1$-valued function of $a$.  Integration against weakly convergent measures gives the first assertion.

For small $r$, local bi-Lipschitz coordinates and a rescaling of the interval also give
\[
 0\le\frac1{2r}\int_{-r}^r\log\frac1{[\psi(x),a]}\,dx
       \le\log\frac1{[z_0,a]}+C,
\]
with $C$ independent of $r,a$.  For the singular part, this reduces to boundedness above of $\log|b|-\tfrac12\int_{-1}^1\log|x-b|\,dx$ for $b\in\C\setminus\{0\}$. Finiteness at $z_0$ permits dominated convergence.  The swept potential and other terms in \eqref{eq:spherical-balayage} are continuous near $z_0$, proving \eqref{eq:arc-mean-at-point}.

Apply Jensen's inequality on any fixed smaller interval $[-r,r]\subset I$:
\[
 \liminf\frac1n\log\int_I f_n(x)\,dx
 \ge-\frac1{2r}\int_{-r}^r H(\psi(x))\,dx.
\]
Let $r\downarrow0$ and use \eqref{eq:arc-mean-at-point} and continuity of $\varphi$.  The assumed upper bound gives the reverse inequality.
\end{proof}

\subsection{Prescribed exceptional poles and weak-convergence rates}\label{app:exceptional-poles}
Here is an elementary construction explaining the scope of the weak pole law. Fix $\beta<\pi$, put $h=1$, and let $0\le k_n\le n$ be integers with $k_n=o(n)$ and $m_n=n-k_n$. Choose a best approximant $b_{m_n}$ and any $k_n$ distinct finite points $a_{j,n}\notin\Sigma_\beta$ avoiding its poles. Write $d_{j,n}=\dist(a_{j,n},\Sigma_\beta)>0$. For $k_n>0$, choose nonzero residues with
\[
 \sum_{j=1}^{k_n}\frac{|c_{j,n}|}{d_{j,n}}\le\rho_\beta^{2n},
 \qquad
 r_n(z)=b_{m_n}(z)+\sum_{j=1}^{k_n}\frac{c_{j,n}}{z-a_{j,n}}.
\]
For $k_n=0$ use $r_n=b_n$. Each new point is a genuine simple pole, because the other terms are holomorphic there; the perturbation is holomorphic at the old poles and cannot cancel their principal parts. The total rational order is at most $m_n+k_n=n$, and
\[
 \eps_n(\beta;1)\le\|e^{-z}-r_n\|_{\Sigma_\beta}
 \le\eps_{m_n}(\beta;1)+\rho_\beta^{2n}.
\]
Since $m_n/n\to1$, Theorem~\ref{thm:root} makes both bounds have $n$th-root limit $\rho_\beta$. Thus prescribed exceptional poles, even with scaled positions approaching the sector boundary, are compatible with root optimality.

There is also no rate depending only on $n$ and $\beta$ for weak convergence over all such sequences. To see this, place the scaled new poles $a_{j,n}/n$ in a fixed small disk disjoint from $\Gamma_\beta\cup\Sigma_\beta$. Take a nonnegative Lipschitz function $f$ on the sphere, zero on $\Gamma_\beta$ and at least $c>0$ on that disk, normalized so that $\|f\|_\infty+\operatorname{Lip}(f)\le1$ for the chordal metric. If $\mu_n$ is the scaled pole measure, the bounded-Lipschitz distance satisfies
\[
 d_{\rm BL}(\mu_n,\nu_\beta)
 :=\sup_{\|g\|_\infty+\operatorname{Lip}(g)\le1}
       \left|\int g\,d(\mu_n-\nu_\beta)\right|
 \ge\int f\,d\mu_n\ge c\,\frac{k_n}{n}.
\]
The fraction $k_n/n$ can tend to zero arbitrarily slowly. Quantitative results may still hold for exact best approximants or with extra assumptions on error loss or exceptional poles.

\section*{Data and code availability}
Two supplementary files accompany this article. Supplementary Material (PDF) contains independent finite-degree bounds, wide-sector asymptotics, scalar-evaluation and numerical-validation details, and matrix-function consequences. The numerical reproduction archive (ZIP) contains fitting and validation scripts, stored barycentric approximants, equilibrium-curve data, and figure-generation scripts. An additional submission-source script computes the new all-pole chordal diagnostic directly from this unchanged archive. The common-denominator and escaping-pole results are proved analytically; no numerical performance claim for those constructions is made.

\section*{Declaration of generative AI and AI-assisted technologies in the manuscript preparation process}
During the preparation of this work, the authors used OpenAI ChatGPT to assist exploratory mathematical derivations, consistency checks, literature comparison, exposition, and code review. This included auxiliary checks associated with the equilibrium construction, boundary localization, local reflection, phase correction, and the common-denominator argument. Anthropic Claude and Perplexity AI were also used for critical reviews of earlier drafts. The mathematical claims are established by the proofs in the article; numerical evidence is identified separately in Section~\ref{sec:num}. The authors reviewed and checked the AI-assisted material, made the final scientific and editorial decisions, and take full responsibility for the content of the article.

\clearpage
\clearpage
\setcounter{section}{0}
\setcounter{equation}{0}
\setcounter{table}{0}
\setcounter{figure}{0}
\setcounter{theorem}{0}
\renewcommand{\thesection}{S\arabic{section}}
\renewcommand{\theequation}{\thesection.\arabic{equation}}
\renewcommand{\thetheorem}{\thesection.\arabic{theorem}}
\renewcommand{\thetable}{S\arabic{table}}
\renewcommand{\thefigure}{S\arabic{figure}}
\renewcommand{\theHsection}{supplement.\arabic{section}}
\renewcommand{\theHtable}{supplement.\arabic{table}}
\renewcommand{\theHfigure}{supplement.\arabic{figure}}

\begin{center}
{\Large\bfseries Supplementary material}\par\medskip
{\large Pole rigidity in root-optimal rational approximation of the exponential on infinite sectors}\par\medskip
Fei Xue and Tianqi Zhang\par
\smallskip
Journal of Mathematical Analysis and Applications\par
\smallskip
Corresponding author: Fei Xue, Clemson University, Clemson, SC 29634, USA\par
\texttt{fxue@clemson.edu}
\end{center}
\medskip
This file collects results and computational details that are not used in the proof of the main article's pole-rigidity theorem. Section S1 gives elementary finite-degree bounds; Section S2 treats the wide-sector limit and numerical evaluation; Section S3 gives matrix-function consequences. All necessary technical proofs for the main theorem remain in the article's appendices.

\section{Finite-degree lower bounds and opening limits}\label{sec:finite-degree}

Throughout,
\[
 \Sigma_\beta=\{re^{i\theta}:r\ge0,\ |\theta|\le\beta/2\}\cup\{\infty\},
 \qquad 0\le\beta<\pi,
\]
and $f_h(z)=e^{-hz}$ for finite $z$, with $f_h(\infty)=0$. The class $\RR_n(\Sigma_\beta)$ consists of reduced rational functions $p/q$ with no sector pole and $\deg p\le\deg q\le n$. Write $\|g\|_E=\sup_E|g|$ and
\[
 \eps_n(\beta;h)=\inf_{r\in\RR_n(\Sigma_\beta)}\|f_h-r\|_{\Sigma_\beta}.
\]
Positive dilation gives $\eps_n(\beta;h)=\eps_n(\beta;1)$. We use the classical Gonchar--Rakhmanov (GR) root factor \cite[\S2.5, p.~329]{supp:GoncharRakhmanov1989} in the form identified in the article:
\begin{equation}\label{eq:root-factor}
 \eps_n(\beta;h)^{1/n}\longrightarrow\rho_\beta
 =\exp\!\left[-\pi\frac{\Kell(1-m_\beta)}{\Kell(m_\beta)}\right],
 \qquad \mathsf P_\eta(1-2m_\beta)=0,
\end{equation}
where $\eta=(\pi-\beta)/(2\pi)$, $1/2<m_\beta<1$, $\mathsf P_\eta$ is the Ferrers function, and
$\Kell(m)=\int_0^{\pi/2}(1-m\sin^2t)^{-1/2}\,dt$ uses the elliptic parameter. Set $v=\rho_0$. The scale parameter $h$ in this document is not a spatial mesh width.

These observations give elementary constraints on the error but play no role in identifying the root-optimal pole measure.

The repeated-pole result in \cite[Theorem~4.5 and Proposition~4.6]{supp:XueZhangRepeated2026} gives the coarse free-pole upper bound $\limsup_{n\to\infty}\eps_n(\beta;1)^{1/n}\le U_\beta$, where $U_\beta=\tan\!\bigl(\pi/[4(2-\beta/\pi)]\bigr)$.  Here $U_\beta$ is the minimum of the repeated-pole rate over the noncoalescing branch analyzed in that work; no stronger repeated-pole global claim is needed below.

\begin{theorem}[Oscillation obstruction]\label{thm:osc}
Let $0<\beta<\pi$.  If $r=p/q$ with $\deg p\le n$, $\deg q\le m$ and $r$ has no pole on $\Sigma_\beta$, then
\[
 \|e^{-z}-r\|_{\Sigma_\beta}
 \ge
 \exp\!\left[-\frac\pi2(n+m+1)\cot\frac\beta2\right].
\]
In particular,
\[
 \eps_n(\beta;1)
 \ge
 \exp\!\left[-\pi\left(n+\frac12\right)\cot\frac\beta2\right].
\]
\end{theorem}

\begin{proof}
Cancel common factors first. The degree bounds do not increase, and the reduced denominator has no zero in the sector. Put $\theta=\beta/2$ and
\[
 X=\frac\pi2(n+m+1)\cot\theta,
 \qquad
 y_j=-\frac\pi2(n+m+1)+j\pi,
 \quad j=0,\ldots,n+m+1.
\]
Then $z_j=X+iy_j$ belongs to $\Sigma_\beta$ and $e^{-z_j}=e^{-X}e^{-iy_0}(-1)^j$.  If the error were strictly smaller than $e^{-X}$, the real numbers
\[
 \Rea\{e^{iy_0}r(X+iy_j)\}
\]
would strictly alternate in sign.  Yet
\[
 \Rea\{e^{iy_0}r(X+iy)\}
 =\frac{\Rea\{e^{iy_0}p(X+iy)\overline{q(X+iy)}\}}
 {|q(X+iy)|^2},
\]
and the numerator is a real polynomial in $y$ of degree at most $n+m$.  The required $n+m+1$ distinct sign changes are impossible.
\end{proof}

For $0<\beta<\pi$, combining Theorem~\ref{thm:osc}, the half-line lower bound $v$, and the repeated-pole construction gives 
\begin{equation}\label{eq:coarse-bounds}
 \max\left\{v,\exp\!\left[-\pi\cot\frac\beta2\right]\right\}
 \le \liminf_n\eps_n(\beta;1)^{1/n}
 \le \limsup_n\eps_n(\beta;1)^{1/n}
 \le U_\beta.
\end{equation}

As an immediate finite-degree consequence, any type-$(n,n)$ approximant with $\|e^{-z}-r\|_{\Sigma_\beta}\le\epsilon<1$ must satisfy $n\ge \pi^{-1}\tan(\beta/2)\log(1/\epsilon)-1/2$.

\begin{proposition}[Opening monotonicity]\label{prop:beta-monotone}
For every fixed $n$, $\eps_n(\beta;1)$ is nondecreasing in $\beta$ and $\lim_{\beta\uparrow\pi}\eps_n(\beta;1)=1$.  Consequently the root factor in \eqref{eq:root-factor} is nondecreasing and satisfies $\rho_0=v$ and $\rho_\beta\uparrow1$ as $\beta\uparrow\pi$.
\end{proposition}

\begin{proof}
The sectors are nested, so the best uniform error is nondecreasing with opening. The zero rational function gives $\eps_n\le1$, whereas Theorem~\ref{thm:osc} gives $\eps_n(\beta;1)\ge\exp[-\pi(n+1/2)\cot(\beta/2)]\to1$.  For the root factor, first take $n$th roots in the oscillation bound and then let $\beta\uparrow\pi$.  Together with the zero-approximant upper bound, this proves the asserted limit without interchanging the two limits.
\end{proof}

\begin{remark}
On the closed right half-plane itself the best finite-degree uniform error is one: along the imaginary axis $e^{-iy}$ winds around the unit circle, while every rational function with finite half-plane sup norm has a finite limit at infinity. Thus $\beta<\pi$ is essential: the limits $n\to\infty$ and $\beta\uparrow\pi$ do not commute.
\end{remark}

\section{Wide-sector asymptotics and computational details}\label{sec:scalar-details}

\subsection{A stable scalar equation}\label{sec:stable-evaluation}
For wide sectors $m_\beta$ is extremely close to one.  Direct floating-point solution in $m$ is therefore ill-conditioned, so set $1-m=e^{-s}$.  The continuation formula \cite[Eq.~15.8.10]{supp:DLMF} rewrites the Ferrers root equation as $G(s)=s+d_0-\sum_{k\ge1}b_ke^{-ks}(s+d_k)=0$, where $\psi$ is the digamma function, $d_k=2\psi(k+1)-\psi(k-\eta)-\psi(k+1+\eta)$ for $k\ge0$, and the positive coefficients satisfy $b_1=\eta(1+\eta)$ and $b_{k+1}/b_k=(k-\eta)(k+1+\eta)/(k+1)^2$.  For $k\ge1$, both $b_k$ and $d_k$ are positive and decreasing.  For the latter, the digamma integral gives positivity, while the recurrence gives $d_{k+1}-d_k=2/(k+1)-1/(k-\eta)-1/(k+1+\eta)<0$.  The tail therefore admits the elementary bound
\[
 0\le\sum_{k>N}b_ke^{-ks}(s+d_k)
 \le
 \frac{b_{N+1}e^{-(N+1)s}(s+d_{N+1})}{1-e^{-s}}.
\]
Together with rigorously enclosed special-function values, this bound allows bisection with outward rounding to certify $s$.  Writing $\operatorname{AGM}$ for the arithmetic--geometric mean, set $A_1=\operatorname{AGM}(1,e^{-s/2})$ and $A_2=\operatorname{AGM}(1,\sqrt{1-e^{-s}})$. Then $-\log\rho_\beta=\pi A_1/A_2$.  The archived code cross-checks complementary-parameter evaluation against these AGM formulas in floating-point arithmetic; Table~\ref{tab:rho} is not interval-certified.

\begin{table}[t]
\caption{The classical free-pole factor and the repeated-pole factor optimized within the noncoalescing branch of \cite{supp:XueZhangRepeated2026}.}
\label{tab:rho}\label{tab:free-vs-repeated}
\centering
\small
\begin{tabular}{c c c@{\quad}c c c}
\toprule
$\beta$ & $\rho_\beta$ & $U_\beta$ & $\beta$ & $\rho_\beta$ & $U_\beta$\\
\midrule
$0^\circ$   & 0.107653919226 & 0.414213562373 &
$120^\circ$ & 0.321876299918 & 0.668178637919\\
$30^\circ$  & 0.133062615880 & 0.456684697904 &
$150^\circ$ & 0.512290800461 & 0.797473388882\\
$60^\circ$  & 0.169677724263 & 0.509525449494 &
$170^\circ$ & 0.775255001842 & 0.920564850211\\
$90^\circ$  & 0.226058759663 & 0.577350269190 &
$179^\circ$ & 0.973160712588 & 0.991359008155\\
\bottomrule
\end{tabular}
\end{table}

\subsection{The wide-sector expansion}
\begin{proposition}[Wide-sector expansion]\label{prop:wide}
Let $\delta=\pi-\beta\downarrow0$.  Then
\begin{equation}\label{eq:wide}
 -\log\rho_\beta
 =\frac\pi2\delta-(\log2)\delta^2+O(\delta^3).
\end{equation}
In particular
\[
 \lim_{\beta\uparrow\pi}
 \frac{-\log\rho_\beta}{\pi-\beta}=\frac\pi2.
\]
\end{proposition}

\begin{proof}
Put $1-m_\beta=e^{-s}$ and $\eta=\delta/(2\pi)$.  Let $\psi$ be the digamma function and $\zeta$ the Riemann zeta function in this proof. The zero-balanced continuation \cite[Eq.~15.8.10]{supp:DLMF}, written in Section~\ref{sec:stable-evaluation}, gives
\[
 s+d_0=\sum_{k\ge1}b_ke^{-ks}(s+d_k),\qquad
 d_0=2\psi(1)-\psi(-\eta)-\psi(1+\eta).
\]
Using $\psi(-\eta)=\psi(1-\eta)+\eta^{-1}$ and the Taylor series \cite[Eq.~5.7.4]{supp:DLMF}, we obtain
\[
 d_0=-\eta^{-1}+2\sum_{j\ge1}\zeta(2j+1)\eta^{2j}
     =-\eta^{-1}+2\zeta(3)\eta^2+2\zeta(5)\eta^4+O(\eta^6).
\]
For small $\eta$, the tail bound in Section~\ref{sec:stable-evaluation}, with $b_1=\eta(1+\eta)$ and bounded $d_1$, brackets the unique root between $-d_0$ and $-d_0+1$.  Uniformly in that interval the sum is $O(e^{-1/\eta})$.  Therefore
\[
 s=-d_0+O(e^{-1/\eta})
   =\eta^{-1}-2\zeta(3)\eta^2+O(\eta^4).
\]
The complementary-parameter expansions in \cite[Sec.~19.12]{supp:DLMF} give
\[
 \Kell(e^{-s})=\frac\pi2+O(e^{-s}),\qquad
 \Kell(1-e^{-s})=\frac{s}{2}+\log4+O(se^{-s}).
\]
Consequently
\[
 -\log\rho_\beta
 =\pi\frac{\Kell(e^{-s})}{\Kell(1-e^{-s})}
 =\frac{\pi^2}{s+4\log2+O(se^{-s})}
 =\pi^2\eta-4\pi^2(\log2)\eta^2+O(\eta^3).
\]
Substituting $\eta=\delta/(2\pi)$ gives \eqref{eq:wide}.
\end{proof}

\subsection{Cross-checks of the equilibrium formulas}
The archived scripts implement two numerical cross-checks separate from rational fitting. Direct Stieltjes quadrature in \texttt{independent\_validation.py} checks the potential at the ODE-computed real crossing for six angles, including $0^\circ$ and $179^\circ$, with maximum residual $1.824\times10^{-12}$. A separate AGM evaluation in \texttt{make\_artifacts.py} checks the elliptic root factor at the eight angles in Table~\ref{tab:rho}, with maximum difference $5.55\times10^{-17}$. Both scripts were prepared by the authors; these are consistency checks by separate implementations, not independent third-party certifications. The complex half-period, mass normalization, and theta--Ferrers identity are proved analytically in the article's Appendices A--B.

\subsection{Fitting and numerical cross-checks}
We set $h=1$ and compute general complex type-$(n,n)$ candidates by adaptive Antoulas--Anderson (AAA) approximation followed by Lawson iteration \cite{supp:NakatsukasaSeteTrefethen2018,supp:NakatsukasaTrefethen2020}.  The accompanying code uses SciPy's AAA implementation and a column-equilibrated, damped Lawson iteration.  The fit uses the coordinate $w=(n-z)/(n+z)$ and samples each boundary ray at $z=nt(1-t)^{-1}e^{\pm i\beta/2}$, with $t=(1-\cos\theta)/2$ and equally spaced $\theta\in[0,\pi]$.  The point $t=1$ represents infinity, with target value zero.  There are $2048$ intervals per ray for $\beta\le150^\circ$ and $8192$ otherwise.  We retain the smallest training error among the initial AAA candidate and $80$ Lawson steps, using weight updates proportional to the square root of the error magnitude. This finite-step procedure provides no guarantee of convergence to a local or global minimizer on the continuous sector.  Neither the S-curve nor predicted poles enter the fit, and no pole--zero pairs are removed.

Errors are evaluated on separate boundary grids with $16384$ and $32768$ intervals per ray, followed by local peak searches.  At $\beta=120^\circ$, $n=24$, we also evaluate the stored double-precision candidate at the refined peaks in $80$-digit arithmetic; the largest peak is checked again at $100$ digits.  The numerical-candidate table in the article uses this higher-precision evaluation for that row.  Every candidate has $n$ finite extracted poles and zero detected sector poles; refined pole residuals and denominator winding counts provide complementary numerical cross-checks.  Source code, full barycentric data, validation records, and figure-generation scripts are included in the supplementary archive.

\subsection{An all-pole chordal support diagnostic}\label{subsec:chordal}
The main article reports a 90th-percentile transversal distance, which deliberately
omits the largest $10\%$ of deviations. To check the entire archived pole set,
write $z_{j,n}=p_{j,n}/n$ and use the normalized chordal metric on the
Riemann sphere,
\[
 [z,w]=\frac{|z-w|}{\sqrt{1+|z|^2}\sqrt{1+|w|^2}},\qquad
 A_n^{\mathrm{ch}}=\frac1n\sum_{j=1}^{n}
       \min_{w\in\Gamma_\beta}[z_{j,n},w].
\]
All stored candidates in the table have exactly $n$ extracted finite poles.
The minimum was evaluated against the piecewise-linear interpolation of the
archived $65\,537$-node equilibrium arc, with one-dimensional minimization on
its nearest segments. Halving the node spacing changes the reported mean by
less than $4\times10^{-10}$ for these cases. The numerical values below are
computed solely from stored candidate poles and stored theoretical arcs; no
AAA--Lawson fit was rerun.
\begin{table}[ht]
\caption{Mean chordal distance of \emph{all} extracted scaled poles to the
predicted S-arc. The lower-degree $179^\circ$ rows expose exceptional transverse
deviations that a fixed percentile can underrepresent.}
\label{tab:all-pole-chordal}
\centering\small
\begin{tabular}{r r r@{\qquad}r r r}
\toprule
$\beta$ & $n$ & $A_n^{\mathrm{ch}}$ & $\beta$ & $n$ & $A_n^{\mathrm{ch}}$\\
\midrule
$90^\circ$ &16&0.009312 &$179^\circ$&64&0.019483\\
$120^\circ$&24&0.005600 &$179^\circ$&96&0.008962\\
$150^\circ$&32&0.002942 &$179^\circ$&128&0.003484\\
$170^\circ$&72&0.000283 &$179^\circ$&240&0.000277\\
$175^\circ$&120&0.000248 & & & \\
\bottomrule
\end{tabular}
\end{table}
For actual root-optimal sequences, weak convergence to $\nu_\beta$ forces
$A_n^{\mathrm{ch}}\to0$, because distance to the compact support is bounded
and continuous on the sphere. The converse is false: this mean checks
transverse proximity but not the equilibrium mass distribution along the arc.
Thus it complements the mass-quantile statistic $D_{\mathrm q}$; the finite-degree
values are diagnostics, not minimax certificates or convergence proofs.
The separate script \texttt{chordal\_diagnostic.py} computes the table from the
existing numerical reproduction archive without changing any fitted candidate.

\section{Matrix-function consequences and degree comparisons}\label{sec:matrix}
Write $W(B)=\{x^*Bx:x^*x=1\}$ for the numerical range of a matrix. In a simple-pole partial-fraction representation, rational evaluation takes the familiar form
\begin{equation}\label{eq:matrix-shifted-solves}
 r(B)b=c_0b+\sum_jc_j(B-p_jI)^{-1}b.
\end{equation}
Thus the poles are computational shifts \cite{supp:Sidje1998,supp:Higham2008}. The following consequence records precisely what a sector-uniform approximation theorem implies; it makes no assumption about the cost of the shifted solves.

\begin{corollary}[Sector-uniform matrix minimax rate and pole law]
\label{cor:matrix-minimax}
For $0\le\beta<\pi$ and $h>0$, define
\begin{equation}\label{eq:matrix-minimax-definition}
 \mathcal E_n^{\mathrm{mat}}(\beta;h)
 :=\inf_{r\in\RR_n(\Sigma_\beta)}
 \sup_{\substack{N\ge1,\ B\in\C^{N\times N}\\W(B)\subset\Sigma_\beta}}
 \|e^{-hB}-r(B)\|_2.
\end{equation}
The same rational function is used for all matrices in the supremum.  Then
\begin{equation}\label{eq:matrix-minimax-sandwich}
 \eps_n(\beta;h)\le\mathcal E_n^{\mathrm{mat}}(\beta;h)
 \le(1+\sqrt2)\eps_n(\beta;h),
 \qquad
 \bigl(\mathcal E_n^{\mathrm{mat}}(\beta;h)\bigr)^{1/n}\to\rho_\beta.
\end{equation}
Moreover, if $r_n\in\RR_n(\Sigma_\beta)$ satisfies
\[
 \left(\sup_{\substack{N\ge1,\ B\in\C^{N\times N},\ W(B)\subset\Sigma_\beta}}
 \|e^{-hB}-r_n(B)\|_2\right)^{1/n}\to\rho_\beta,
\]
then its reduced pole measures satisfy the pole law of the main article: $n^{-1}\sum_p\delta_{hp/n}\weakto\nu_\beta$, where $\nu_\beta$ is its arc-supported equilibrium measure.
\end{corollary}

\begin{proof}
For each admissible $r$, scalar matrices $B=[z]$, with finite $z\in\Sigma_\beta$, give the lower bound by $\|f_h-r\|_{\Sigma_\beta}$; the value at infinity is its limit along the positive ray.  The Crouzeix--Palencia estimate \cite{supp:CrouzeixPalencia2017} gives the upper bound by $(1+\sqrt2)\|f_h-r\|_{\Sigma_\beta}$.  For each matrix the function $f_h-r$ is holomorphic on a neighborhood of its compact numerical range because $r$ has no sector poles.  Taking the supremum over matrices and then the infimum over $r$ proves the sandwich; \eqref{eq:root-factor} gives the root limit.  For the final assertion,
\[
 \eps_n(\beta;h)\le\|f_h-r_n\|_{\Sigma_\beta}
 \le\sup_{\substack{N\ge1,\ B\in\C^{N\times N}\\W(B)\subset\Sigma_\beta}}
 \|e^{-hB}-r_n(B)\|_2.
\]
Thus $r_n$ is scalar root-optimal, and the main article's pole-rigidity theorem applies.
\end{proof}

For a fixed opening, let $n_{\rm free}(\epsilon)$ be the least degree for which the scalar best error is at most $\epsilon$. Inverting the root limit gives the approximation-degree benchmark
\begin{equation}\label{eq:matrix-work-count}
 n_{\rm free}(\epsilon)
 =\frac{\log(1/\epsilon)}{-\log\rho_\beta}(1+o(1)),
 \qquad \epsilon\downarrow0.
\end{equation}
The same asymptotic degree holds for the sector-uniform matrix error because the sandwich constant is independent of degree. This counts rational degrees, not wall-clock work: a repeated pole can reuse a factorization, while distinct shifts have different conditioning and solution costs.

For $\delta=\pi-\beta\downarrow0$, $-\log U_\beta=\delta/2+O(\delta^2)$.  Together with Proposition~\ref{prop:wide}, this makes the ratio of the stated repeated-pole branch's degree estimate to the free-pole degree estimate tend to $\pi$, with accuracy tending to zero first at fixed angle.  This is not a joint angle--degree estimate or a runtime ratio.

\paragraph{Files and scope.}
The numerical reproduction archive is a separate supplementary file. The statements above are either proved here or are explicit consequences of the main article and cited literature; none is a missing step in the main theorem. All empirical claims refer to the stored candidates and validation records, not to global minimax or interval certificates.

\end{document}